\documentclass[bj,preprint]{imsart}

\usepackage[autonum]{tchdr}
\usepackage{natbib}
\usepackage{xcolor}
\newcommand{\prokd}[2]{d_{\mathrm{P}}(#1,#2)}
\startlocaldefs

\endlocaldefs

\allowdisplaybreaks

\begin{document}

\begin{frontmatter}
\title{Local Exchangeability}
\runtitle{Local Exchangeability}

\begin{aug}
\author[A]{\fnms{Trevor} \snm{Campbell}\ead[label=e1,mark]{trevor@stat.ubc.ca}},
\author[A]{\fnms{Saifuddin} \snm{Syed}\ead[label=e2,mark]{saif.syed@stat.ubc.ca}},
\author[B]{\fnms{Chiao-Yu} \snm{Yang}\ead[label=e3,mark]{chiaoyu@berkeley.edu}},
\author[B]{\fnms{Michael I.} \snm{Jordan}\ead[label=e4,mark]{jordan@eecs.berkeley.edu}},
\and
\author[C]{\fnms{Tamara} \snm{Broderick}\ead[label=e5,mark]{tbroderick@mit.edu}}

\address[A]{Department of Statistics,
University of British Columbia, Vancouver, Canada.
\printead{e1,e2}}

\address[B]{Department of Electrical Engineering and Computer Science,
University of California Berkeley, Berkeley, USA.
\printead{e3,e4}}

\address[C]{Laboratory for Information and Decision Systems,
Massachusetts Institute of Technology, Cambridge, USA.
\printead{e5}}

\runauthor{T. Campbell et al.}

\end{aug}

\begin{abstract}
Exchangeability---in which the distribution of an infinite sequence
is invariant to reorderings of its elements---implies the existence of 
a simple conditional independence structure that may be leveraged
in the design of statistical models and inference procedures.
In this work, we study a relaxation of exchangeability in which
this invariance need not hold precisely.
We introduce the notion of \emph{local exchangeability}---where swapping
data associated with nearby covariates causes a bounded change in the distribution.
We prove that locally exchangeable processes correspond to independent observations 
from an underlying measure-valued stochastic process. 
Using this main probabilistic result, we show that the \emph{local empirical measure}
of a finite collection of observations provides an approximation of the 
underlying measure-valued process and Bayesian posterior predictive distributions.
The paper concludes with applications of the main theoretical results
to a model from Bayesian nonparametrics and covariate-dependent permutation tests.
\end{abstract}


\end{frontmatter}
\section{Introduction}
Let $X = X_1, X_2, \dots$ be an infinite sequence of random elements in
a standard Borel space $(\mcX, \Sigma)$.
The sequence is said to be \emph{exchangeable}
if for any finite permutation $\pi$ of $\nats$,
\[
X_1, X_2, \dots \eqd X_{\pi(1)}, X_{\pi(2)}, \dots. \label{eq:exch}
\]
At first sight this assumption appears innocent; intuitively, it suggests only
that the order in which observations appear provides no information about those
or future observations.  But despite its apparent innocence, exchangeability
has a powerful implication.  In particular, the well-known 
\emph{de Finetti's theorem} \citep[e.g.][Theorem 11.10]{Kallenberg02} states that an infinite
sequence is exchangeable if and only if it is mixture of \iid sequences, i.e.,
there exists a unique random probability measure $G$ on $\mcX$ such that
\[
\Pr\left(X\in \cdot \given G\right) \eqas G^{\infty}, \label{eq:deF}
\]
where $G^\infty$ is the countable infinite product measure constructed from
$G$.  Thus, exchangeability provides a strong justification for the Bayesian
approach to modeling \citep{Jordan10}, and guarantees a latent conditional
independence structure of $X$ useful in the design of computationally efficient
inference algorithms.  Exchangeability is also the basis of well-known
nonparametric permutation testing procedures
(\citealp{Pitman37a,Pitman37b,Pitman37c}; \citealp[Ch.~3]{Fisher66}; \citealp{Ernst04}; \citealp[Ch.~15]{Lehmann05}).

However, although exchangeability may be a useful idealization in modeling and
analysis, many data come with covariates that preclude an honest
belief in its validity. For example, given a corpus of documents tagged by
publication date, one might reasonably expect the data to exhibit
a time-dependence that is incompatible with exchangeability. 
Nevertheless, one might still expect the distribution not to change
too much if we permuted documents published only one day apart;
i.e., observations with similar covariates are intuitively ``nearly
exchangeable.''
In this work, we investigate how to codify this intuition.

One option is to use a kind of \emph{partial exchangeability} \citep{deFinetti38,Lauritzen74,Diaconis78,Camerlenghi19} 
in which the distribution is invariant to permutations within equivalence classes. Formally, we
endow each observation $X_n$ with a covariate $t_n$ from a set $\mcT$,
and assert that the sequence distribution is invariant only to reordering observations with
equivalent covariate values.
Under this assumption as well as the availability of infinitely many observations at each covariate value, 
we have a similar representation of $X$ as a mixture of independent
sequences given random probability measures $(G_{t})_{t\in \mcT}$,
\[
\Pr\left(X \in \cdot \given (G_{t})_{t\in \mcT}\right) \eqas \prod_{n=1}^\infty G_{t_n}. \label{eq:deFPartial}
\]
The random probability measures $(G_t)_{t\in \mcT}$ can have an arbitrary dependence on one another;
partially exchangeable sequences encompass those that are exchangeable
(where the covariate does not matter), decoupled (where subsequences for each different
covariate value are mutually independent), and the full range of models in between. 
In particular, partial exchangeability does not enforce the desideratum that observations with 
nearby covariates should have a similar law, and is too weak to be useful for restricting the class of 
underlying mixing measures for the data.

In this work, we introduce a new notion of \emph{local exchangeability}---lying
between partial and exact exchangeability---in which swapping data associated
with nearby covariates causes a bounded change in total variation distance.  We
begin by studying probabilistic properties of locally exchangeable processes in
\cref{sec:defn,sec:definetti}. The main result from this section is in the spirit of de
Finetti's theorem: we prove that locally exchangeable processes correspond to
independent observations from a unique underlying smooth measure-valued
stochastic process.  To the best of our knowledge, this representation theorem
is the first to arise from an approximate probabilistic symmetry.  Further, the
existence of such an underlying process not only shows that de Finetti's
theorem is robust to perturbations away from exact exchangeability, justifying
the Bayesian analysis of real data, but also imposes a useful constraint on the
space of models one should consider when dealing with data that one suspects
follows a locally exchangeable random process.
Next in \cref{sec:localempiricalmeasures}, we
use this result to show that the \emph{local empirical measure} 
of a finite collection of observations can be used to provide 
an approximation of the
underlying measure-valued process,
Bayesian predictive posterior distributions, and 
the premetric that governs local exchangeability.
These results rely heavily on the intuition that 
locally exchangeable observations from nearby covariates 
behave essentially like exchangeable observations.
Finally, in \cref{sec:examples}, we provide 
example applications in two statistical models exhibiting
local exchangeability---Gaussian processes \citep{Rasmussen06} 
and dependent 
Dirichlet processes \citep{MacEachern99,MacEachern00}---as well as grouped permutation 
tests in the presence of covariates.
The paper
concludes with a discussion of directions for future work.  Proofs of all
results are provided in the appendix.

\subsection{Related work}
Beyond de Finetti's original result for infinite binary sequences \citep{deFinetti31} and its
extensions to more general range spaces \citep{deFinetti37,Hewitt55} and 
finite sequences \citep{Diaconis77,Diaconis80}---see \citet{Aldous85} for 
an in-depth introduction---correspondences between probabilistic invariances and conditional
latent structure (known as \emph{representation theorems}) have been studied extensively.
Notions of exchangeability and corresponding latent conditional structure
now exist for a wide variety of probabilistic models,
such as arrays \citep{Aldous81,Hoover79,Austin14,Jung18},
Markov processes \citep{Diaconis80b},
networks \citep{Caron17,Veitch15,Borgs16,Crane16,Cai16,Janson17},
combinatorial structures \citep{Kingman78,Pitman95,Broderick13,Campbell18,Crane16b},
random measures \citep{Kallenberg90},
and more \citep{Diaconis88,Kallenberg05,Orbanz15}.
Furthermore, weaker notions of exchangeability such as conditionally identical distributions \citep{Berti04,Kallenberg88} have been developed. All past work on probabilistic invariance and its consequences has pertained 
to exact invariance.

\section{Local exchangeability}\label{sec:localexch}
\subsection{Definition}\label{sec:defn}
Let $X = (X_t)_{t\in\mcT}$ be a stochastic process on an index (or covariate) set $\mcT$ taking
values in a standard Borel space $(\mcX, \Sigma)$.  To encode distance between
covariates, we endow the set $\mcT$ with a \emph{premetric} $d: \mcT \times
\mcT \to [0,1]$ satisfying $d(t, t') = d(t', t)$ and $d(t, t) = 0$ for $t, t'
\in \mcT$.  
We will formalize local exchangeability based on the finite dimensional
projections of $X$.  For any subset $T \subset \mcT$ and injection $\pi : T \to \mcT$, 
let 
$X_T$ and $X_{\pi, T}$ denote stochastic processes on index set $T$ such that
\[
\forall t \in T, \qquad (X_T)_t \defined X_t \qquad (X_{\pi, T})_t \defined X_{\pi(t)}. \label{eq:XpiTdefn}
\]
In other words, $X_T$ is the restriction of $X$ to index set $T$, while $X_{\pi, T}$ 
is the restriction to $T$ under the mapping $\pi$. \cref{defn:localexch} captures 
the notion that observations with similar
covariates should be close to exchangeable, i.e., the total variation between
$X_T$ and $X_{\pi, T}$ is small as long as the distances between $t$ and $\pi(t)$
are small for all $t\in T$.
\bnumdefn\label{defn:localexch}
The process $X$ is
\emph{locally exchangeable} with respect to a premetric $d$ if
for any finite subset $T\subset \mcT$
and injection $\pi: T \to \mcT$,
\[
	d_{\mathrm{TV}}(X_T, X_{\pi, T}) \leq \sum_{t\in T}d(t, \pi(t)). \label{eq:localexch}
\]
\enumdefn
\cref{defn:localexch} generalizes both exchangeability and partial exchangeability among equivalence classes. In 
particular, the zero premetric where $d(t,t') = 0$ identically yields classical
exchangeability, while the premetric $d(t, t') = 1 - \ind[t \sim t']$ for equivalence relation $\sim$
yields partial exchangeability. Further, any process is locally exchangeable with respect to
the discrete premetric $d(t,t') = 1- \ind[t = t']$; in order to say something of
value about a process $X$, it must satisfy \cref{eq:localexch} for a tighter
premetric.

To quantify differences in distributions, \cref{defn:localexch} employs the total variation distance,
which for random elements $Y,Z$ in a
measurable space $(\mcY, \Xi)$ is defined as
\[
d_{\mathrm{TV}}(Y,Z) \defined \sup_{A\in\Xi} \left| \Pr(Y\in A) - \Pr(Z \in A)\right|.
\]
The choice of total variation distance (as opposed to other metrics and divergences, see e.g.~\citep{Gibbs02})
is motivated by its symmetry and generality. 
We make $d$ a premetric---as opposed to a (pseudo)metric, say---as 
 the triangle inequality and positive definiteness are unused in the theory below.
Further we use a premetric with range $[0, 1]$ because total variation always lies in
this range, and so any valid bound in \cref{eq:localexch} for 
a premetric $d : \mcT \times \mcT \to \reals_+$ can be improved by replacing $d$ with
$\min(d, 1)$. 
And although \cref{defn:localexch} imposes a total variation bound only for all finite sets 
of covariates, it is equivalent to do so for all countable sets of covariates,
as shown in \cref{prop:finiteinfiniteequiv}.  
\bnprop\label{prop:finiteinfiniteequiv}
If $X$ is locally exchangeable with respect to $d$, then 
for any countable subset $T \subset \mcT$
and injection $\pi: T \to \mcT$,
\[
	d_{\mathrm{TV}}(X_{T}, X_{\pi, T}) \leq \sum_{t\in T}d(t, \pi(t)). \label{eq:inflocalexch}
\]
\enprop

\bnexa\label{eg:linreg} A simple example of local exchangeability that we will return to throughout the paper is
the process of observable measurements $X$ 
from a Bayesian linear regression model on $\mcT = \reals$ with a quadratic trend,
\[
\theta \dist \distNorm(0, 1), \qquad \forall t\in\reals, \quad X_t \distind \distNorm(\theta t^2, 1). \label{eq:bayeslinreg}
\]
By \cref{lem:prodmsrbd},
since the $X_t$ are independent conditioned on $\theta$, 
\[
d_{\mathrm{TV}}(X_T, X_{\pi,T}) &\leq \sum_{t\in T} \EE\left[ d_{\mathrm{TV}}(\distNorm(\theta t^2, 1), \distNorm(\theta \pi(t)^2, 1)) \right].
\]
We bound the terms in the sum
using the Lipschitz continuity of the standard normal CDF $\Phi$,
\[
\EE\left[d_{\mathrm{TV}}( \distNorm(\theta t^2, 1), \distNorm(\theta\pi(t)^2, 1 ))\right]
&= \EE\left[\Phi\left( \frac{|\theta t^2 - \theta\pi(t)^2|}{2} \right) - 
\Phi\left( -\frac{|\theta t^2 - \theta \pi(t)^2|}{2} \right)\right]\\
&\leq \frac{\EE|\theta||t^2 - \pi(t)^2|}{\sqrt{2\pi}}\leq \frac{|t^2 - \pi(t)^2|}{\sqrt{2\pi}}.
\]
Therefore the process $X$ in the Bayesian linear regression model \cref{eq:bayeslinreg} 
is locally exchangeable with respect to the premetric $d(t,t') = \min(|t^2 - t'^2|/\sqrt{2\pi}, 1)$.
Note that
we are free to take $\min(\cdot, 1)$ because the total variation is bounded above by 1.
This example illustrates why we opt for the generality of a premetric; here, observations at points $t$ and $-t$
are exactly exchangeable since $d(t,-t) = 0$, which does not generally hold for a metric, and
$|t^2 - t'^2|$ does not satisfy the triangle inequality.
Also note that the marginal distribution of $X_T$ is a multivariate Gaussian with off-diagonal covariance terms 
$\EE\left[X_t X_{t'}\right] \propto t^2 t'^2$, which varies with $t, t'$;
multivariate Gaussians with exchangeable components must have constant off-diagonal covariance terms. 
Therefore this example also shows that there exist processes that
are locally exchangeable but not exchangeable.
\enexa

\subsection{de Finetti representation} \label{sec:definetti}
In the previous example, we used the fact that the variables $X_t$ were conditionally
independent given a latent random variable $\theta$ to demonstrate their local exchangeability.
A natural question to ask is whether \emph{all} locally exchangeable processes
exhibit a similar structure. 
\cref{thm:dfrep} answers this question in the affirmative, by providing a de Finetti-like representation
of locally exchangeable processes similar to \cref{eq:deF} and \cref{eq:deFPartial}.
This representation guarantees the existence
of a simple conditional structure that can be leveraged in the design of statistical inference
procedures, and justifies
a Bayesian approach when dealing with covariate-dependent data. 
We first require a weak assumption on the space $\mcT$.
\bnumdefn[Infinitely-separable space]\label{assum:infsep} A premetric space $(d, \mcT)$ is \emph{infinitely separable} if 
there exists a countable subset
$\mfT \subseteq \mcT$ such that for all $t\in \mcT$, there exists a Cauchy sequence 
$(t_n)_{n\in\nats}$ in $\mfT$ such 
that $t_n \to t$ and $|\{t_n : n\in\nats\}| = \infty$.
\enumdefn
When $d$ is a metric, infinite separability is equivalent to $\mcT$ being separable with no isolated points. 
When $d$ is a pseudometric,
it is equivalent to the existence of a countable dense subset $\mfT \subseteq \mcT$
such that for all $t\in \mcT$ and $\epsilon > 0$, $|\{t' \in \mfT : d(t,t') < \epsilon\}| = \infty$.
In general, infinite separability ensures that 
there are infinitely many elements to swap ``nearby'' each covariate value of interest $t\in\mcT$. 
This assumption precludes the situation where observations 
satisfy finite exchangeability \citep{Diaconis77,Diaconis80} but not infinite exchangeability.

\cref{thm:dfrep}
shows that under infinite separability, the desired de Finetti-like representation indeed does exist. 
In particular, we show that there is a unique probability measure-valued process $G$ that renders $X$
conditionally independent, and that $G$ satisfies a
continuity property with the same ``smoothness'' as the observed process.
For the precise statement of the result in \cref{thm:dfrep}, 
recall that a \emph{modification} of a stochastic process $G$
on $\mcT$ is any other process $G'$ on $\mcT$ such that $\forall t\in\mcT, \, \Pr\left(G_t = G'_t\right) = 1$.
\bnthm\label{thm:dfrep}
Suppose $(d, \mcT)$ is infinitely separable. 
Then the process $X$ is locally exchangeable with respect to $d$ if and only if
there exists a random measure-valued stochastic process $G=(G_t)_{t\in\mcT}$
(unique up to modification)
such that for any finite subset of covariates $T\subset \mcT$ and $t, t' \in \mcT$,
\[
\Pr\left(X_T\in\!\cdot\! \given G\right) &\eqas \prod_{t\in T}G_t,  &
\sup_{A}\EE\left| G_t(A) - G_{t'}(A) \right|&\leq d(t, t').\label{eq:aldeFrep}
\]
\enthm
For example, given $\mcT = \nats$ and the zero premetric $d(t,t') = 0$, one recovers 
the de Finetti representation of exchangeable sequences; the smoothness 
condition asserts that $G_t$ must be constant for all $t\in\mcT$ as expected. 
Similarly, suppose we are given an equivalence relation $\sim$ 
on $\nats$ where each equivalence class has infinite cardinality. Then setting $\mcT = \nats$
and $d(t,t') = 1 - \ind[t \sim t']$ recovers the de Finetti representation of partially
exchangeable sequences under permutation within equivalence classes; here the smoothness condition
asserts that $G_t$ must be constant within each equivalence class, but allows for general dependence
between $G_t$ across the equivalence classes. Thus, in the same way
that \cref{defn:localexch} generalizes (partial) exchangeability, \cref{thm:dfrep} generalizes
the de Finetti representation theorem. 

Note that we still obtain the ``if'' direction of \cref{thm:dfrep} without 
imposing the infinite separability assumption on $(d,\mcT)$. 
In particular, if we are given a process $G$
satisfying \cref{eq:aldeFrep}, 
then the process $X$ is locally exchangeable with respect to 
both 
\[
d_c(t,t')&\defined \sup_A \EE\left|G_t(A) - G_{t'}(A)\right|,
&
\text{and}&
&
d_{sc}(t,t')&\defined \EE\left[ d_{\mathrm{TV}}(G_t, G_{t'}) \right]. \label{eq:canonicalmetric}
\]
We refer to $d_c$ as the \emph{canonical premetric} and
$d_{sc}$ as the \emph{strong canonical premetric}. 
Note that $X$ is locally exchangeable
with respect any premetric $d$ satisfying $d \geq d_c$, and in particular, $d_{sc} \geq d_c$.
Given a particular $G$, one can 
use \cref{lem:prodmsrbd}  to derive 
an upper bound on these two premetrics (as demonstrated in \cref{eg:linreg}),
which then provides insight into the extent to which data $X$ 
generated from $G$ are exchangeable.
Note that $(d_c, \mcT)$ and $(d_{sc}, \mcT)$ may or may
not be infinitely separable, depending on the characteristics of the process $G$.

\bexa[continued]
In the linear regression example, the underlying measure-valued process is 
the collection of normal distributions 
\[
G_t = \distNorm(\theta t^2, 1), \qquad t \in \mcT.
\]
\cref{thm:dfrep} guarantees that this process is unique up to modification.
In this case, the randomness in $G$ is entirely due to the latent variable $\theta \dist \distNorm(0,1)$;
in general $G$ need not be determined by a finite-dimensional quantity.
We can also verify that $G$ satisfies the required smoothness condition
with respect to $d$, although it is not surprising in this case given that we originally
derived the premetric using the same technique:
\[
\sup_A \EE\left| G_t(A) - G_{t'}(A)\right| &\leq \EE\, d_{\mathrm{TV}}(G_t, G_{t'}) \leq \min\left(\frac{1}{\sqrt{2\pi}}|t^2 - t'^2|, 1\right) = d(t,t').
\]
\eexa

\subsection{Local empirical measure process}\label{sec:localempiricalmeasures}
The de Finetti result in \cref{thm:dfrep} guarantees the existence of a unique underlying
measure-valued process $G$, but does not provide any direct insight 
into the distribution of $G$ or whether it is identifiable
given only (countably many) measurements of the process $X$.
In the classical setting of an exchangeable sequence $X_1, X_2, \dots$, the empirical measure 
$\whG_N = \frac{1}{N}\sum_{n=1}^N \delta_{X_n}$ of 
a finite collection of observations $(X_n)_{n=1}^N$ 
serves this purpose, as it converges weakly to $G$ almost surely \citep{Varadarajan58}, i.e.,
\[
\prokd{\whG_N}{G}\convas 0,\qquad N\to\infty,\label{eq:exchempmsr}
\]
where $d_\mathrm{P}$ denotes the L\'evy-Prokhorov metric.
In the setting of local exchangeability more generally, however, the usual empirical measure
does not provide a result similar to \cref{eq:exchempmsr}. If we are interested in understanding the distribution of
$G_\tau$ for some $\tau\in\mcT$, and we collect measurements $(X_t)_{t\in T}$
of $X$ at a finite set of covariates $T\subset \mcT$, the presence of far-away
covariates in $T$ from $\tau$ can result in a non-vanishing bias in the
empirical measure. 
To address this issue, for each $\tau \in \mcT$, let $t_i(\tau)$, $i=1,\dots, |T|$ be an ordering
of the set $T$ such that the values $d_i(\tau) = d(t_i(\tau), \tau)$ 
are ordered from smallest to largest.
Then define
\[
M_\tau &= \max\left\{M \in [|T|] : \frac{1}{M}\left(1+\sum_{m=1}^M2 d_m(\tau)\right) > 2d_M(\tau)\right\} ,
& 
\mu_\tau &=\frac{1}{M_\tau}\sum_{m=1}^{M_\tau} d_m(\tau). 
\]
We construct
the \emph{local empirical measure process} $(\whG_\tau)_{\tau\in\mcT}$ via
\[
\whG_\tau &= \sum_{t\in T}\xi_t(\tau) \delta_{X_t}, \qquad \xi_t(\tau) 
= \max\left\{0 ,\frac{1}{M_\tau}+2(\mu_\tau - d(t,\tau))\right\}.  \label{eq:whg}
\]
The local empirical measure process $\whG$ serves as an approximation of the
measure-valued process $G$ underlying the locally exchangeable process $X$.
Note that
$\sum_{t\in T}\max\{0 ,\frac{1}{M_\tau} + 2(\mu_\tau - d(t,\tau))\} = 1$, so $\whG_\tau$ is 
a probability measure for each $\tau\in\mcT$.
Further note that $(\whG)_{\tau\in\mcT}$ is measurable with respect to $(X_t)_{t\in T}$.
Intuitively, $\whG$ includes only those observations at covariates sufficiently 
close to the point of interest $\tau\in\mcT$ such that the decrease in variance
associated with adding another observation outweighs the potential increase in bias.
The value $M_\tau$ represents
how many observations are included in the local empirical measure at that location,
and $\mu_\tau$ represents the average distance of their covariates to $\tau$.

Our goal now is to provide a weak convergence result for the local empirical measure process $\whG$ in the limit of
many observations, similar to that of \cref{eq:exchempmsr}.  
As a key step towards that goal, \cref{thm:empiricalmeasurebds} provides bounds on both the
expected squared estimation error (\cref{eq:whgsqerrorbd}) as well as error
tail probabilities (\cref{eq:whgprbd}) when using the local empirical measure process
$\whG_\tau$ in place of $G_\tau$ or $\Pr\left(X_\tau \in \cdot \given X_T\right)$, for all $\tau\in\mcT$.  
Each bound in \cref{thm:empiricalmeasurebds} has
two terms: the first is related to the variance incurred by estimation via
independent sampling, and the second is related to the bias incurred by using
observations from $t\neq \tau$.
Note that \cref{thm:empiricalmeasurebds} 
quantifies the approximation error using the metric
\[
\|\nu-\eta\|_{\mcA} = \sum_{i=1}^\infty c_i \left|\nu(A_i) - \eta(A_i)\right|, \quad \nu,\eta\text{ probability measures}, 
\]
where $\mcA = \left\{c_i, A_i\right\}_{i=1}^\infty$, 
$A_i$ are measurable subsets of $\mcX$, $c_i \geq 0$, and $\sum_i c_i = 1$.
We work with $\|\cdot\|_\mcA$ rather than standard metrics because it simplifies the analysis substantially.
Although the properties of $\|\cdot\|_\mcA$ depend on the choice of 
$\mcA$ in general, there exists a choice such that $\|\cdot\|_\mcA \to 0$ implies weak
convergence (see \cref{lem:normAimpliesweakconv} in the appendix),
and the bounds below in \cref{thm:empiricalmeasurebds} are valid for any choice of $\mcA$, as indicated by the supremum.
We will use the metric $\|\cdot\|_\mcA$ and the results in \cref{thm:empiricalmeasurebds}
as a stepping stone to obtain weak convergence in \cref{cor:empiricalmeasureasymp} below.

\bnthm\label{thm:empiricalmeasurebds}
Let $(d, \mcT)$ be infinitely separable and $X$ be locally exchangeable with respect to $d$. 
Then 
\[
\forall \tau \in \mcT, \quad \sup_{\mcA} \EE\left[\|\whG_{\tau} - G_{\tau}\|^2_\mcA\right] &\leq 
\frac{1}{4M_\tau} + \mu_\tau, \label{eq:whgsqerrorbd}
\]
and for all $\delta > 0$, $\tau\in\mcT$,
\[
\sup_{\mcA} \Pr\left(\|\whG_\tau - G_\tau\|_\mcA > \delta + \sqrt{2\mu_\tau + 1/M_\tau}\right)&\leq 
\exp\left( \frac{-\delta^2}{2\left(2\mu_\tau+1/M_\tau\right)}\right)
+\frac{2\mu_\tau}{\delta + \sqrt{1/M_\tau}}.
\label{eq:whgprbd}
\]
Furthermore, the same bounds in \cref{eq:whgsqerrorbd,eq:whgprbd} apply 
when $G_\tau$ is replaced with $\Pr\left(X_\tau \in \cdot \given X_T\right)$.
\enthm
When all of the covariates in the observed set $T$ are close to $\tau$,
the bounds in \cref{thm:empiricalmeasurebds} provide essentially the same 
guarantees as one would expect for exchangeable random variables. In 
particular, suppose for all $t\in T$, $d(t, \tau) \lesssim \exp(-|T|)$, and so $\xi_t(\tau) \approx 1/|T|$. 
In this situation the bounds above reduce to
\[
\sup_\mcA \EE\left[\|\whG_{\tau} - G_{\tau}\|^2_\mcA\right] = O(|T|^{-1}), \quad 
\sup_\mcA \Pr\left(\|\whG_{\tau} - G_{\tau}\|_\mcA > \delta + |T|^{-1/2}\right) &= O\left(e^{-|T|\delta^2}\right).
\]

\cref{cor:empiricalmeasureasymp} uses the results in \cref{thm:empiricalmeasurebds} to 
obtain a weak convergence result for $\whG_\tau$ similar to \cref{eq:exchempmsr}.
In particular, if we collect measurements of $X$ from a 
sequence of sets that concentrate around $\tau$---for example,
$T_n = \{t_i\}_{i=1}^n$ such that there exists a subsequence $t_{i_k} \to \tau$---then the
local empirical measure $\whG_\tau$ converges weakly to both $G_\tau$ 
and the Bayesian posterior predictive distribution in probability.
Recall that $d_\mathrm{P}$ denotes the L\'evy-Prokhorov metric. 

\bncor\label{cor:empiricalmeasureasymp}
Fix $\tau \in \mcT$. Suppose we make observations at a sequence of finite sets $T_n \subset \mcT$, $n\in\nats$
of covariates such that for all $\epsilon > 0$, $\left|\{t \in T_n : d(t,\tau) \leq \epsilon\}\right| \to \infty$.
Then
 \[
 \prokd{\whG_\tau}{G_\tau} \convp 0
\quad\text{and}\quad
 \prokd{\whG_\tau}{\Pr\left(X_\tau \in \cdot \given X_{T_n}\right)}\convp 0,
\qquad n\to\infty.
\]
\encor
A byproduct of \cref{cor:empiricalmeasureasymp} is that one 
can characterize the distribution of $G_\tau$ by analyzing
the distribution of $X_\tau$ conditioned on
$X_{T_n}$ for a sequence of sets of 
covariates $T_n$ that concentrate around $\tau$,
i.e., $|T_n|\to\infty$ and 
$\max\{d(t, \tau) : t\in T_n\} \to 0$ as $n\to\infty$.
Note that it is not required to know the premetric $d$ governing local exchangeability
in order to identify $G$ using this technique;
one can instead construct the set of covariates $T_n$ such that $\max\{\ell(t, \tau) : t\in T_n\} \to 0$
for any premetric $\ell : \mcT \times \mcT \to [0,1]$
that \emph{dominates} $d$ in the sense that 
for any two sequences of covariates $t_n, t'_n$, $n\in\nats$, 
\[
\ell(t_n, t'_n) \to 0 \implies d(t_n, t'_n) \to 0, \quad n \to \infty.\label{eq:ldominated}
\]
The requirement in \cref{eq:ldominated} is typically not stringent; 
it states only that when covariates get close under $\ell$,
they must also get close under $d$, with no other stipulation about relative rates, bounds, etc.
In the following linear regression example, we
will use the usual metric $\ell(t,t') = |t - t'|$ on $\reals$.
\bexa[continued] We return to the linear regression
example to show how the distribution of $G_\tau$ 
can be recovered from the process $X$ via \cref{cor:empiricalmeasureasymp}.
The joint density of $X_T, X_\tau$ is
\[
p(x_\tau, x_T) &\propto \exp\left(-\frac{1}{2}x_\tau^2 -\frac{1}{2}\sum_{t\in T} x_t^2 + \frac{1}{2}\frac{\left(x_\tau\tau^2 + \sum_{t\in T} x_t t^2\right)^2}{1 + \tau^4+ \sum_{t\in T} t^4} \right).
\]
Therefore the conditional distribution of $X_\tau$ given $X_T$ is given by
\[
X_\tau \dist \distNorm\left(\frac{\tau^2\sum_{t\in T}X_t t^2}{1+\sum_{t\in T}t^4}, \frac{1+\tau^4+\sum_{t\in T}t^4}{1+\sum_{t\in T}t^4}\right), 
\]
If we then consider a sequence of sets $T_n$ of 
covariates that grows in size and concentrates quickly
around $\tau$---e.g.,
$T_n = \{\tau + i\exp(-n) : i=1, \dots, n\}$---we
find that the conditional distribution of $X_\tau$ given $X_T$ converges to
\[
X_\tau \dist \distNorm\left(Y, 1\right), 
\,\,\text{where}\,\,
Y \dist \distNorm\left(0, \tau^4\right).
\]
By setting $\theta = Y\tau^{-2}$, we recover the fact 
that $X_\tau$ is generated from $G_\tau = \distNorm(\theta\tau^2, 1)$,
$\theta \dist\distNorm(0, 1)$, i.e., 
the marginal of the original Bayesian linear regression model. 
Note that
one can repeat essentially the same analysis
for multiple covariates $\tau_1, \dots, \tau_K$ to recover finite marginal distributions. 
For example, if we consider the bivariate distribution of $G_{\tau_1}, G_{\tau_2}$,
we find that $X_{\tau_1}, X_{\tau_2}$ are generated independently from
\[
G_{\tau_1} = \distNorm(\theta\tau_1^2,1) \quad G_{\tau_2} = \distNorm(\theta\tau_2^2, 1), \quad  \theta\dist\distNorm(0,1).
\]
The analysis from the example in \cref{sec:defn} can then be used to bound the 
strong canonical 
premetric $d_{sc}(t,t') = d_{\mathrm{TV}}(G_t, G_{t'}) \leq \min\left(|t-t'|/\sqrt{2\pi}, 1\right)$.
Thus, given only the process $X$, we have 
identified a premetric $d$ under which $X$ is locally exchangeable
as well as the measure-valued process $G$.
\eexa

\subsection{Regularity}\label{sec:regularity}
The smoothness property of $G$ in \cref{eq:aldeFrep} may seem unsatisfying at a first glance; it 
bounds the absolute difference in the underlying mixing measure process at 
nearby locations only \emph{in expectation}, leaving room for the possibility of sample discontinuities
in $G_t$ as a function of $t$.  However, there are many probabilistic models 
that, intuitively, generate observations that should be considered locally exchangeable but which have discontinuous 
latent mixing measures. For example, some dynamic nonparametric mixture models \citep{Lin10,Chen13} have components that
are created and destroyed over time, causing discrete jumps in the mixing measure. As long as the jumps happen 
at diffuse random times, the probability of a jump occurring between two times decreases
as the difference in time decreases, and the observations may still be locally exchangeable. 
However, intuitively, if there is a fixed location $t_0$ with a nonzero probability of a discrete jump in the mixing
measure process, the observations $X$ cannot be locally exchangeable. \cref{cor:diffusejumps} provides the precise statement.
\bncor\label{cor:diffusejumps}
Suppose $(d, \mcT)$ is infinitely separable and $X$ 
is locally exchangeable with respect to $d$. 
Then for all $A\in\Sigma$, $t_0\in\mcT$, and $\epsilon>0$,
\[
\lim_{\eta \to 0}\sup_{t \,:\, d(t,t_0)\leq \eta}\Pr\left(|G_t(A) - G_{t_0}(A)| > \epsilon\right) = 0.
\] 
\encor
That being said, it is worth examining whether different guarantees on properties of the underlying measure process $G$
result as a consequence of different properties of the premetric $d$. \cref{cor:stationarity} answers this question in the affirmative
for processes on $\mcT=\reals$;
in particular, the faster the decay of $d(t, t')$ relative to $|t-t'|$ as $t\to t'$, the stronger the guarantees on the behavior
of the mixing measure $G$. Note that while this result is presented for covariate space $\reals$,
the result can be extended to processes on $\reals\times \nats$ and more general 
separable spaces \cite[Theorems 2.8, 2.9, 4.5]{Pothoff09}.
\bnthm\label{cor:stationarity}
Let $\mcT = \reals$, $\gamma \geq0$, and $X$ be locally exchangeable 
with respect to a premetric $d$ satisfying $d(t, t') = O(|t-t'|^{1+\gamma})$ as $|t-t'| \to 0$.
Then:
\benum
\item ($\gamma > 1$): $X$ is exchangeable and $G$ is a constant process.
\item ($0<\gamma\leq 1$): $X$ is stationary and for any $A\in\Sigma$ and $\alpha \in (0, \gamma)$,  $(G_t(A))_{t\in\reals}$ is weak-sense stationary with an $\alpha$-H\"older continuous modification.
\item ($\gamma = 0$): $G$ may have no continuous modification.
\eenum
\enthm
\brmk
A rough converse of the first point holds: $X$ exchangeable implies constant $G$, and $d(t,t') = 0$ is trivially $O(|t-t'|^{1+\gamma})$ for $\gamma > 1$.
But a similar claim for the second point is not true in general: $X$ stationary and locally exchangeable does not necessarily imply that $d(t,t')=O(|t-t'|^{1+\gamma})$ for $0<\gamma\leq 1$.
 For a counterexample, consider a square wave shifted by a uniform random variable,
i.e., the process $X_t = \mathrm{sign}\left(\sin(2\pi(t-U))\right)$ 
for $U\dist\distUnif[0, 1]$. Here
$X_t$ is stationary and locally exchangeable with $d(t,t') = \min(|t-t'|,1)$, 
but $|t-t'| \neq O(|t-t'|^{1+\gamma})$ for any
$\gamma > 0$ as $|t-t'|\to 0$.
\ermk

\subsection{Approximate conditional independence}
In the classical setting of exchangeable sequences $X_1, X_2, \dots$,
the empirical measure $\whG = \frac{1}{N}\sum_{n=1}^N \delta_{X_n}$ satisfies the following property:
for all bounded measurable functions $h:\mathcal{X}^N\to\mathbb{R}$,
\[
\EE\l[h(X_1,\dots,X_N)|\whG, G\r]=\EE\l[h(X_1,\dots,X_N)|\whG \r].\label{eq:suff}
\]
Thus $G$ and $(X_1, \dots, X_N)$ are conditionally independent given $\whG$.
In other words, the fact that $(X_1, \dots, X_N)$ corresponds to covariate values $(1, \dots, N)$ 
provides no additional information about $G$ beyond $\whG$ itself.

In the setting of local exchangeability, the question of how important the covariate values are in inferring the measure-valued process $G$
is relevant in practice: we do not often get to observe the true covariate values $\{t_1, \dots, t_N\} = T\subset \mcT$,
but rather we observe discretized versions that are grouped into ``bins.'' For example, if $X_T$ corresponds to observed document data
 with timestamps $T$, we may know those timestamps up to only a certain precision (e.g.~days, months, years).
This section shows that a ``binned'' version of the empirical measure $\whG$ provides an approximate conditional independence similar
to \cref{eq:suff},
where the error of approximation decays smoothly by an amount corresponding to the uncertainty in covariate values.

Formally, suppose we partition our covariate space $\mathcal{T}$ into disjoint bins 
$\{\mcT_k\}_{k=1}^\infty$, where each bin has observations $T_k=\mathcal{T}_k\cap T$. We may use a finite partition 
by setting all but finitely many $\mcT_k$ to the empty set.
Although we know the number of points in each bin (i.e., the cardinality of $T_k$), we will encode our lack of knowledge of their positions as randomness: $T_k\sim \mu_k$, where $\mu_k$ is a probability distribution capturing
our belief of how the unobserved covariates are generated within each bin. 
Following the intuition from the classical de Finetti's theorem, we define the binned
empirical measures $\wtG_{k}=\sum_{t\in T_k}\delta_{X_t}$, $\wtG\defined (\wtG_1,\wtG_2,\dots)$, and
let $\mcG$ denote the subgroup of permutations $\pi:T\to T$ that permute observations only within each bin, i.e., such that $\forall k\in\nats$, $\pi(T_k)=T_k$.
Note that  $|\mcG|=\prod_{k=1}^\infty |T_k|! < \infty$ since there are only finitely many observations in total. Unlike classical exchangeability, 
 $\wtG$ does not provide exact conditional independence of $X_T$ and $G$; but \cref{prop:approxsuff} guarantees that 
 it provides a form of approximate conditional independence, with error that depends on $(\mu_k)_{k=1}^\infty$.
\bnthm\label{prop:approxsuff}
Suppose $(d, \mcT)$ is infinitely separable. If $X$ is locally exchangeable with respect to $d$, and $h:\mathcal{X}^T\to\mathbb{R}$ is a bounded measurable function,
\[
\EE\left|\EE\left[h(X_T)\given \wtG,G\right]-\EE\left[h(X_T)\given \wtG\right]\right|\leq 4\|h\|_\infty\EE\left[\sum_{t\in T}d(t,\pi(t))\right],\label{eq:approxsuff}
\]
where $\pi\dist\distUnif\left(\mcG\right)$ and $T_k \distind \mu_k$.
\enthm
\brmk
Note that the expectation on the right hand side  
averages over the randomness both in the uncertain covariates $T$ and the permutation $\pi$.
\ermk
If $X$ is exchangeable within each bin $\mathcal{T}_k$, \cref{prop:approxsuff}
states that $X_T$ and $G$ are conditionally independent given $\wtG$, as desired. Further, the deviance from
independence is controlled by the deviance from exchangeability within each bin.
In particular, 
\[
\EE\left[\sum_{t\in T} d(t,\pi(t))\right]
&\leq \sum_{k=1}^\infty |T_k|\diam{\mathcal{T}_k} \leq |T|\sup_k\{\diam{\mathcal{T}_k}\},\label{eq:diambds}
\]
where $\diam{\mcT_k} \defined \sup_{t,t' \in \mcT_k}d(t,t')$.
Both bounds in \cref{eq:diambds} are independent of $\mu_k$; thus 
the result holds even if we are unwilling to express our uncertainty in the binned covariates via a
distribution.

\section{Examples}\label{sec:examples} 

In this section, we provide example applications of the theory in
\cref{sec:localexch}.  First, we use a case study of Gaussian processes to show
how one can use posterior predictive distributions to analyze the local
exchangeability of a process. In particular, we show how to derive 
the underlying measure process $G$, as well as
an appropriate premetric $d$ governing local 
exchangeability, using only finite marginals of the process $X$.
Second, we use a case study of dependent Dirichlet processes to show that 
one can use local empirical measures as a surrogate
for otherwise intractable posterior predictive distributions in discrete Bayesian nonparametric models.
See the appendix for other examples of Bayesian nonparametric models
exhibiting local exchangeability---e.g., kernel beta process feature models
\citep{Hjort90,Ren11} and dynamic topic models \citep{Blei06,Wang08}, among
others. Finally, we demonstrate a usage of local exchangeability as a tool to
analyze the inflation of type-I error in matched permutation tests involving covariates.

\subsection{Obtaining the underlying measure-valued process and premetric}
We will first provide
an example of how one can use the Bayesian posterior predictive distributions
of a locally exchangeable process $X$ to derive the distribution of the underlying
measure-valued process $G$ as well as the premetric of local exchangeability $d$.
This example applies the same strategy as in 
the running example from \cref{sec:localempiricalmeasures}, albeit in 
a more sophisticated nonparametric model.

Consider a Gaussian process  $X\dist\distGP(m, \kappa)$ on $\mcT = \reals^d$ with 
continuous mean function $m : \reals^d \to \reals$,
and covariance function $\kappa(x,y) = \sigma^2(x) \ind[x=y] + k(x,y)$
for continuous nonnegative $\sigma^2 : \reals^d \to \reals_+$ and
continuous symmetric positive-definite kernel $k : \reals^d\times\reals^d \to \reals_+$.
Define a set of $k$ unique covariate values 
$\tau_1, \dots, \tau_k \in \mcT$, and consider the Euclidean metric on $\mcT$. 
For each $n\in\nats$ and $i=1,\dots, k$, 
let $T_{in}$ be a finite subset of covariates
such that $|T_{in}| = n$ and 
$\max\{\|\tau_i-t\| : t\in T_{in}\} = o(1/n)$.
Direct analysis of the conditional density 
yields that as $n\to\infty$, the conditional distribution 
of $X_{\tau_1}, \dots, X_{\tau_k}$ given $X_{T_{1n}}, \dots, X_{T_{kn}}$
converges to
\[
(X_{\tau_1}, \dots, X_{\tau_k})
\dist \distNorm\left(\left(Y_1, \dots, Y_k\right), \diag\left(\sigma^2(\tau_1),\dots, \sigma^2(\tau_k)\right)\right), \label{eq:gp1}
\]
where
\[
\left(Y_1, \dots, Y_k\right)\dist\distNorm\left(\left(m(\tau_1),\dots,m(\tau_k)\right), K\right), \qquad K_{ij} = k(\tau_i,\tau_j). \label{eq:gp2}
\]
\cref{eq:gp1,eq:gp2} demonstrate 
that $X$ is conditionally independently drawn from the process $G$ where
\[
\forall \tau\in\mcT, \quad G_\tau = \distNorm(Y_\tau, \sigma^2(\tau)) \qquad Y \dist\distGP(m, k).
\]
We now derive the strong canonical premetric of local exchangeability.
In this setting, 
\[
d_{sc}(t,t') = \EE\left[d_{\mathrm{TV}}(G_t, G_{t'})\right] = \EE\left[d_{\mathrm{TV}}( \distNorm(Y_t, \sigma^2(t)), \distNorm(Y_{t'}, \sigma^2(t'))) \right].\label{eq:dctvgp}
\]
By \citet[Theorem 1.3]{Devroye20},
\[
d_{\mathrm{TV}}( \distNorm(Y_t, \sigma^2(t)), \distNorm(Y_{t'}, \sigma^2(t'))) 
&\leq \frac{3|\sigma^2(t)-\sigma^2(t')|}{2\max\{\sigma^2(t),\sigma^2(t')\}} + \frac{|Y_t - Y_{t'}|}{2\max\{\sigma(t),\sigma(t')\}}.
\]
Applying Jensen's inequality $\EE|Y_t - Y_{t'}| \leq \sqrt{\EE(Y_t - Y_{t'})^2}$,
then evaluating the expectation and using the bounds
$|\sigma^2(t) - \sigma^2(t')| \leq 2\max\{\sigma(t),\sigma(t')\}|\sigma(t)-\sigma(t')|$,
and $\sqrt{x^2+y^2} \leq x+y$ yields
\[
d_{sc}(t,t')
&\leq \min\left(1, \frac{6|\sigma(t)-\sigma(t')| + |m(t)-m(t')| + \sqrt{k(t,t) + k(t',t') - 2k(t,t')}}{2\max\{\sigma(t),\sigma(t')\}}\right). \label{eq:gpboundref}
\]
In the usual setting with zero mean $m(t) = 0$, constant noise variance
$\sigma(t) = \sigma$ for some $\sigma>0$, and stationary 
kernel $k(t,t') = r(\|t-t'\|)$ for some $r: \reals_+\to\reals_+$, 
\cref{eq:gpboundref} reduces to
\[
d_{sc}(t,t') 
&\leq \min\left(1, \frac{\sqrt{r(0) - r(\|t-t'\|)}}{\sigma}\right).
\]
This example demonstrates that Gaussian processes are locally
exchangeable in the presence of measurement noise, i.e.~where $\sigma(t) > 0$.
However, note that $\sigma(t) > 0$ is not strictly necessary for local exchangeability; 
to obtain a
necessary and sufficient characterization of local exchangeability in Gaussian
processes, we could instead analyze the canonical metric $d_c$ per \cref{thm:dfrep}.

\subsection{Approximate predictive distributions in discrete Bayesian nonparametrics}
Next, we demonstrate that the local empirical measure can serve as a useful surrogate
for otherwise intractable posterior predictive
distributions in discrete Bayesian nonparametric models.
The Dirichlet process \citep{Ferguson73}
is a popular prior for the weights and component parameters in nonparametric
mixture models. Draws from a Dirichlet process are discrete probability measures, 
\[
G = \sum_{k=1}^\infty w_k \delta_{\theta_k}, 
\]
where $(w_k)_{k=1}^\infty$ are weights satisfying $w_k \geq 0$, $\sum_k w_k = 1$,
and $(\theta_k)_{k=1}^\infty$ are component parameters,
each with distribution given by \citep{Sethuraman94}
\[
\theta_k\distiid H, \qquad v_k\distiid \distBeta(1, \alpha), \qquad w_k = v_k\prod_{i=1}^{k-1}(1-v_i), \qquad k&\in\nats,
\]
for some distribution $H$ and concentration parameter $\alpha > 0$. 
Given draws $X_n \distiid G$, the posterior predictive distribution 
of $X_{N+1}$ given the first $N$ draws $X_{1}, \dots, X_{N}$ is
\[
X_{N+1} \dist \frac{\alpha}{\alpha+N} H + \frac{1}{\alpha+N}\sum_{n=1}^N \delta_{X_n} = \frac{\alpha}{\alpha+N}H + \frac{N}{\alpha+N}\whG. \label{eq:dppostpred}
\]
The fact that one can marginalize the (infinitely many) weights and parameters to arrive at \cref{eq:dppostpred} 
is critical in tractable computational inference for models involving the Dirichlet process \citep{Neal00}.

When the observations come with additional covariate information,
the \emph{dependent} Dirichlet process mixture model \citep{MacEachern99,MacEachern00} may be used instead. 
There are many instantiations of the dependent Dirichlet process; for simplicity we consider a model
where the weights are a function of a covariate but the component parameters are 
constant across covariate values, i.e.,
\[
X_{x,n}&\distind \sum_{k=1}^\infty w_{x,k} \delta_{\theta_{k}}, \qquad n\in\nats, \,\,x\in\reals,
\]
where $w_{x,k} = v_{x,k}\prod_{i=1}^{k-1}(1-v_{x,i})$, and 
the stick variables $v_{x,k}$ are now \iid stochastic processes on $\reals$. 
The marginal distributions of $v_{x,k}$ at 
$x\in\reals$ are designed to be $\distBeta(1,\alpha)$
so that the dependent Dirichlet process is marginally a Dirichlet 
process for each covariate value. 
But even for simple stochastic processes $v_{x,k}$, the posterior predictive distribution
is not tractable to obtain in closed-form.
However, we can note that the process $X$ is locally exchangeable with strong canonical premetric
\[
d_{sc}(t,t') &= \EE\left[d_{\mathrm{TV}}\left(\sum_{k=1}^\infty w_{x,k} \delta_{\theta_k}, \sum_{k=1}^\infty w_{x',k} \delta_{\theta_k}\right)\right] = \frac{1}{2}\sum_{k=1}^\infty \EE\left|w_{x,k} - w_{x',k}\right|,
\]
where $t = (x,n)$ and $t'=(x',n')$. Since $w_{x,k}$ is a product of independent variables, \cref{lem:prodbd} yields
\[
d_{sc}(t,t') &\leq 
\frac{1}{2}\EE\left[\left|v_{x,1}-v_{x',1}\right|\right]\sum_{k=1}^\infty 
\left(\left(\frac{\alpha}{\alpha+1}\right)^{k-1} +\frac{k-1}{1+\alpha}\left(\frac{\alpha}{\alpha+1}\right)^{k-2}\right).
\]
The infinite sum converges to some $0<C<\infty$, and so
\[
d_{sc}(t,t') &\leq \min\left(1, C\EE\left|v_{x,1}-v_{x',1}\right|\right).
\]
Therefore, as long as the stochastic process $v_{x,1}$ is smooth enough,
and we condition on $X_T$, where $T$ concentrates closely around 
$\tau \in \mcT$, the posterior predictive distribution of 
$X_\tau$ given $X_T$ 
is approximately equal to the local empirical measure $\whG_\tau$, by \cref{thm:empiricalmeasurebds}; 
the latter has a tractable closed-form expression.

\subsection{Type-I error inflation in grouped permutation tests}\label{sec:localrandomization}
One of the key applications of exchangeability in statistical data analysis is
in the design of nonparametric \emph{permutation tests} with exact type-I
error bounds (\citealp{Pitman37a,Pitman37b,Pitman37c}; \citealp[Ch.~3]{Fisher66}).  
In the notation of this work, we are given observations of a stochastic process $X$ at a finite set 
of covariates $T\subset \mcT$,
a subgroup of $\mcG$ permutations $\pi: T \to T$,  
and a test statistic $S : \mcX^T \to \reals$.
The null hypothesis is that $X_T$ is exchangeable; so
we set a desired threshold $\alpha \in[0, 1]$, and reject the null 
with type-I error at most $\alpha$ if
\[
\frac{1}{|\mcG|}\sum_{\pi \in \mcG}\ind\left[S(X_{T}) \leq S(X_{\pi,T})\right] \leq \alpha,
\]
where $X_{\pi,T}$ is defined as in \cref{eq:XpiTdefn}.
This setup is commonly used in observational studies 
with a control group and treatment group,
where $\mcG$ consists of permutations that swap \emph{matched pairs} 
of elements in the control and treatment 
groups. However, a typical problem is that elements in the two groups are not exactly comparable
due to the presence of covariates. In this case, a standard approach is to 
construct $\mcG$ to permute only those elements 
with similar covariates from the control and treatment 
groups, under some metric $d$ \citep{Cochran65,Rubin73,Rubin73b,Rosenbaum89,Rosenbaum02,Lu04,Greevy04,Hansen04,Hansen06,Baiocchi10,Lu11}. 
Local exchangeability provides a general way to 
analyze the type-I error of these methods; \cref{prop:exactrandomizationtest} shows 
 that for a locally exchangeable process, 
the type-I error $\alpha$ may potentially be increased by the average
distance between pairs of covariates permuted by $\pi\in\mcG$.
\cref{eq:localrandomizationbd} also incidentally provides a rigorous justification
for past work that formulates the construction of $\mcG$ as the minimization
of this penalty (e.g., \citet{Rosenbaum89}).
\bnprop\label{prop:exactrandomizationtest}
Let $X$ be locally exchangeable with respect to $d$. For $\alpha \in [0,1]$,
\[
\Pr\left( \frac{1}{|\mcG|}\sum_{\pi \in \mcG}\ind\left[S(X_T) \leq S(X_{\pi, T})\right] \leq \alpha \right) 
\leq \alpha + \frac{1}{|\mcG|}\sum_{\pi\in\mcG}\sum_{t\in T}d(t, \pi(t)).\label{eq:localrandomizationbd} 
\]
\enprop

\section{Discussion}
The major question posed in this paper is what we can do with data when we do
not believe that they are exchangeable, but are willing to believe that they
are \emph{nearly} exchangeable. This paper answers the question with a
relaxed notion of \emph{local} exchangeability in which swapping data
associated with nearby covariates causes a bounded change in total variation
distance. We have demonstrated that classical results for exchangeable processes are
``robust to the real world;'' indeed, locally exchangeable processes have a de
Finetti representation that may be 
leveraged in the design of statistical models and inference procedures.
Finally, many popular covariate-dependent statistical models---which
violate the assumptions of exchangeability---satisfy local exchangeability,
extending the reach of exchangeability-based analyses to these models.

One limitation of local exchangeability is the infinite separability assumption.  
There are  applications in which the covariate space $\mcT$ has isolated
points that violate this condition, e.g.,
discrete time series where the covariate space is $\mcT = \nats$
endowed with the Euclidean metric. However, if $X$ can be extended to a process
on $\mcS \supseteq \mcT$ such that $(d, \mcS)$ is infinitely separable and
 $(X_s)_{s\in\mcS}$ is locally exchangeable with respect to $d$,
then the theoretical results from this work hold for the marginal process $(X_t)_{t\in\mcT}$.
Another limitation is that the total variation bound in the definition of local
exchangeability is quite weak, which has downstream consequences for the
tightness of the error bounds in \cref{sec:localempiricalmeasures}.
Further study on alternate definitions of
local exchangeability is warranted to strengthen these guarantees.

As a final note, it is also possible that an analogue of the theory of finite
exchangeability \citep{Diaconis80} holds in the local setting; but it is not
yet clear whether this is indeed true or what form it would take.
It would also be of interest to investigate more general notions of
local exchangeability under group actions, e.g., permutations that
 preserve some statistic of the data, which have been used 
in past work on randomization testing in the presence of covariates \citep{Rosenbaum84}.

\section*{Acknowledgements}
The authors thank Jonathan Huggins for illuminating discussions.
T.~Campbell is supported by a National Sciences and Engineering Research Council of Canada (NSERC) Discovery 
Grant and Discovery Launch Supplement. T.~Broderick is supported in part by an NSF CAREER Award, an ARO YIP Award,
ONR, and a Sloan Research Fellowship.

\appendix

\section{Proofs}\label{sec:proofs}
\bprfof{\cref{prop:finiteinfiniteequiv}}
Choose some ordering of the countable set $T = (t_1, t_2, \dots)$. We note that
$(X_{t_n})_{n=1}^\infty$ and
$(X_{\pi(t_n)})_{n=1}^\infty$ are $\mcX^\infty$-valued random
variables that are measurable with respect to $\Sigma^\infty$, 
which is generated by the algebra of cylinder sets
of the form $U \times \mcX^\infty$ for $U\in\Sigma^N$. 
Therefore, 
\[
d_{\mathrm{TV}}(X_T, X_{\pi, T}) &= 
\sup_{N\in\nats, U\in\Sigma^N} 
\left| \Pr\left(\left(X_{t_n}\right)_{n=1}^N \in U\right) - \Pr\left(\left(X_{\pi(t_n)}\right)_{n=1}^N \in U\right)\right|,
\]
where we have replaced $\Sigma^\infty$ with its generator
by the fact that for any algebra of sets $\mcA$,
$\epsilon > 0$, $B \in \sigma(\mcA)$, and probability measures
$\mu, \nu$ on $\sigma(\mcA)$, there exists an $A\in\mcA$
such that $\frac{1}{2}(\mu+\nu)(B\triangle A) < \epsilon$.
So by the definition of local exchangeability for finite sets of covariates,
\[
d_{\mathrm{TV}}(X_T, X_{\pi, T}) 
&\leq \sup_{N\in\nats} \sum_{n=1}^N d(t_n, \pi(t_n)) \leq \sum_{t \in T} d(t, \pi(t)).
\]
\eprfof

\bprfof{\cref{thm:dfrep}}
We start with the reverse direction.
Define the two product measures 
$G_T = \prod_{t\in T}G_t$ and $G_{\pi, T} = \prod_{t \in T} G_{\pi(t)}$.
Then since $\Pr(X_T\in A) = \EE\left[G_T(A)\right]$
and $\Pr(X_{\pi, T}\in A) = \EE\left[G_{\pi, T}(A)\right]$,
by Jensen's inequality,
\[
\sup_{A}\left|\Pr(X_T\in A) - \Pr(X_{\pi, T} \in A)\right|
&= \sup_{A} \left|\EE\left[G_T(A) - G_{\pi, T}(A)\right]\right|\\
&\leq \sup_{A} \EE\left|G_T(A) - G_{\pi, T}(A)\right|\\
&= \sup_A \EE\left| \int \ind[x\in A] \prod_{t\in T}\dee G_t(x_t) - \int \ind[x\in A] \prod_{t\in T}\dee G_{\pi(t)}(x_t)  \right|.
\]
Finally, the proof technique of \citet[Lem.~2.1]{Sendler75}
and the smoothness of $G$
yields the conclusion,
\[
&\leq \sum_{t \in T} \sup_A \EE\left| G_t(A) - G_{\pi(t)}(A)\right|\\
&\leq \sum_{t\in T}d(t, \pi(t)).
\]
For the forward direction, suppose $X$ is locally exchangeable.
Let $(t_n)_{n=1}^\infty$ be any ordering of the countable set $\mfT$ from \cref{assum:infsep},
and let $\mcF$ be the tail $\sigma$-algebra of $(X_{t_n})_{n=1}^\infty$.
We will show that for any two covariates $r, s \in \mcT\setminus \mfT$, $r\neq s$,
$X_r$ and $X_s$
are conditionally independent given $\mcF$. The argument extends via standard methods
to $r,s$ that may be elements of $\mfT$, and then to any finite subset of $\mcT$.

By infinite separability (\cref{assum:infsep}), there exists a subsequence $i_1<i_2< \dots$ of indices such
that $t_{i_n}$ is Cauchy and converges to $s$. By taking another subsequence we can assume without loss of generality for all $N\in\nats$, $i_N > N$ and $d(s,t_{i_{N}})+\sum_{n=N}^\infty d(t_{i_n}, t_{i_{n+1}}) < 1/N$.
Let $\pi_N$ be the mapping that takes
$s\to t_{i_N}$, $t_{i_n} \to t_{i_{n+1}}$ for all $n \geq N$,
and leaves all other $t\in\mcT$ fixed.
Then denote
$Y_N = (X_s, X_{t_N}, X_{t_{N+1}}, \dots)$,
and let $Z_N$ be the sequence with covariates mapped under $\pi_N$.
By reverse martingale convergence, for any bounded measurable $\phi : \mcX \to \reals$,
\[
\EE\left[\phi(X_r) \given Y_N\right] &\convas \EE\left[\phi(X_r) \given X_s, \mcF \right]\quad\text{and}\quad
\EE\left[\phi(X_r) \given Z_N\right] \convas \EE\left[\phi(X_r) \given \mcF \right]
\]
as $N\to\infty$.
Next, by local exchangeability and \cref{prop:finiteinfiniteequiv},
\[
d_{\mathrm{TV}}( (X_r, Y_N), (X_r, Z_N) ) < \frac{1}{N},
\]
and 
by \cref{lem:conditionaldists}(2), we have that
the Wasserstein distance between 
$\EE\left[\phi(X_r) \given Y_N\right]$ and $\EE\left[\phi(X_r) \given Z_N\right]$
converges to 0 as $N\to\infty$. Together, the Wasserstein distance bound and reverse
martingale above yield
\[
\EE\left[\phi(X_r) \given X_s, \mcF \right] \eqd\EE\left[\phi(X_r) \given \mcF \right].
\]
By \citet[Lemma 3.4]{Aldous85},
\[
\EE\left[\phi(X_r) \given X_s, \mcF \right] \eqas \EE\left[\phi(X_r) \given \mcF \right],
\]
and thus $X_r$ and $X_s$ are conditionally independent given $\mcF$.
As mentioned earlier this argument extends to any finite subset $T$ of covariates,
by considering subsequences of $(t_n)_{n=1}^\infty$ converging to each $t\in T$.
Since $X$ takes values in a standard Borel space, there is a random measure $G_t$
for each $t\in\mcT$ for which $G_t(A) \eqas \EE\left[\ind[X_t\in A] \given \mcF\right]$
 \citep[e.g.][Theorem 6.3]{Kallenberg02}. The collection of these random measures
forms the desired stochastic process $G=(G_t)_{t\in\mcT}$.

Next, we develop the smoothness property of $G$. By both reverse and forward martingale convergence, we have that
\[
\MoveEqLeft{\sup_{A}\EE\left|G_t(A) - G_{t'}(A)\right|}\\
 &=\sup_{A}\EE\left|\lim_{n\to\infty}\lim_{m\to\infty} \EE\left[ \ind\left[X_t \in A\right] - \ind\left[X_{t'}\in A\right] \given X_{t_{n:n+m}}\right]\right|, 
\]
Using dominated convergence to move the limits out of the expectation, local exchangeability
to bound the total variation between $(X_t, X_{t_{n:n+m}})$ and $(X_{t'}, X_{t_{n:n+m}})$, and \cref{lem:conditionaldists}(1),
\[
\sup_{A}\EE\left|G_t(A) - G_{t'}(A)\right| &\leq d(t, t'). \label{eq:glocalexch}
\]

Finally, we show that $G$ is approximated by empirical averages of the observations $X$;
this property will be used below to show that $G$ is unique up to modification. 
Consider any $A\in\Sigma$ and any sequence $(t'_n)_{n=1}^\infty$ converging to $s\in \mcT$ such that $d(t'_n, s) \leq 2^{-n}$ for each $n\in\nats$. Define $S_{s,N} = \frac{1}{N}\sum_{n=1}^N \ind[X_{t'_{n}} \in A]$.
Then
\[
\MoveEqLeft{\Pr\left( |S_{s,N} - G_s(A)| > \epsilon\right) }\\
&=\EE\left[\Pr\left( \left|S_{s,N} - \frac{1}{N}\sum_{n=1}^N G_{t'_n}(A) + \frac{1}{N}\sum_{n=1}^N G_{t'_n}(A) - G_s(A)\right| > \epsilon \given \mcF\right)\right]\\
&\leq\EE\left[\Pr\left( \left|S_{s,N} - \frac{1}{N}\sum_{n=1}^N G_{t'_n}(A)\right| + \left|\frac{1}{N}\sum_{n=1}^N G_{t'_n}(A) - G_s(A)\right| > \epsilon \given \mcF\right)\right].
\]
Noting that the right term is $\mcF$-measurable and applying Hoeffding's inequality to the left, 
\[
\Pr\left( |S_{s,N} - G_s(A)| > \epsilon\right) \leq
\EE\left[ 2e^{-2N\left(\max\left\{0, \epsilon - \left|\frac{1}{N}\sum_{n=1}^N G_{t'_n}(A) - G_s(A)\right|\right\}\right)^2}  \right].
\]
Splitting the above expectation across two events---one where the measures satisfy 
\[
\left|\frac{1}{N}\sum_{n=1}^N G_{t'_n}(A) - G_s(A)\right| > \nicefrac{\epsilon}{2}
\]
and the other its complement---yields
\[
\Pr\left( |S_{s,N} - G_s(A)| > \epsilon\right) \leq
\Pr\left( \left|\frac{1}{N}\sum_{n=1}^N G_{t'_n}(A) - G_s(A)\right| > \frac{\epsilon}{2}\right)
+ 
2e^{-N\epsilon^2/2}  .
\]
Applying Markov's inequality, the triangle inequality, and \cref{eq:glocalexch},
\[
\Pr\left( |S_{s,N} - G_s(A)| > \epsilon\right) &\leq
\frac{2}{\epsilon N} \sum_{n=1}^{N}\EE\left|G_{t'_n}(A) - G_s(A)\right| + 2e^{-N\epsilon^2/2}\\
&\leq\frac{2}{\epsilon N} \sum_{n=1}^{N} 2^{-n} + 2e^{-N\epsilon^2/2}
\to 0, \quad \text{as}\quad N\to\infty.
\]
Thus, $S_{s,N} \convp G_s(A)$. We now show that $G$ is unique.
Suppose there is another measure process $G'$ that satisfies \cref{eq:aldeFrep}, from which $X$ is generated 
conditionally independently given some $\sigma$-algebra $\mcF'$.
By repeating the steps above, one can show that $S_{s,N} \convp G'_s(A)$.
Therefore, 
\[
\forall A \in \Sigma, \quad \Pr\left(G_s(A) = G'_s(A)\right) = 1.
\]
Since $(\mcX, \Sigma)$ is a standard Borel space, $\Sigma = \sigma(\mcA)$ for some countable
algebra of sets $\mcA$ \cite[Prop.~3.1, 3.3]{Preston08}. By noting 
that the countable intersection of unit-measure sets is also
unit-measure,
\[
\Pr\left(\forall A \in \mcA, G_s(A) = G'_s(A)\right) = 1.
\]
Finally by Carath\'eodory's 
extension theorem \citep[Theorem 2.5]{Kallenberg02}, the probability measures $G_s$ and $G'_s$ are almost surely equal.
 The extension of this argument 
to any finite subset of covariates $T \subset \mcT$
is straightforward, implying that $(G_t)_{t\in\mcT}$ is uniquely determined
up to modification.
\eprfof

\bprfof{\cref{thm:empiricalmeasurebds}}
First, since $c_i \geq 0$ and $\sum_i c_i = 1$, by Jensen's inequality,
\[
\EE\left[\|\whG_{\tau} - G_{\tau}\|^2_\mcA\right] &= \EE\left[\left(\sum_{i=1}^\infty c_i |\whG_\tau(A_i) - G_\tau(A_i)|\right)^2\right]\\
&\leq  \sum_{i=1}^\infty c_i\EE\left[(\whG_{\tau}(A_i) - G_{\tau}(A_i))^2\right].
\]
We will focus on a single term in the sum for some $A\in\Sigma$ and drop the $i$ subscript, as the bound for all terms will be identical.
Adding and subtracting $\sum_{t\in T}\xi_t(\tau)G_t(A)$,
\[
\EE\left[(\whG_{\tau}(A) - G_{\tau}(A))^2\right]
&=\EE\left[\left(\whG_{\tau}(A) - \sum_{t\in T}\xi_t(\tau)G_t(A) + \sum_{t\in T}\xi_t(\tau)G_t(A) - G_{\tau}(A)\right)^2\right].
\]
Since $X$ is locally exchangeable, by \cref{thm:dfrep}, it is conditionally independently drawn from $G$.
Therefore $\EE\left[\whG_\tau(A) \given G\right] = \sum_{t\in T}\xi_t(\tau)G_t(A)$.
Hence we can use the tower property and expand the square to find that
\[
\EE\left[(\whG_{\tau}(A) - G_{\tau}(A))^2\right]
&=
\EE\left[\left(\whG_{\tau}(A) - \sum_{t\in T}\xi_t(\tau)G_t(A)\right)^2\right] 
+ \EE\left[\left(\sum_{t\in T}\xi_t(\tau)G_t(A) - G_\tau(A)\right)^2\right].
\]
The first term can be bounded by using the same conditional independence property again---in particular,
that $\EE\left[ \ind\left[X_t \in A\right] \given G\right] = G_t(A)$---followed by Popoviciu's inequality:
\[
\EE\left[\left(\whG_{\tau}(A) - \sum_{t\in T}\xi_t(\tau)G_t(A)\right)^2\right] &=
\EE\left[\left(\sum_{t\in T}\xi_t(\tau)\left(\ind[X_t \in A] - G_t(A)\right)\right)^2\right]\\
&=\sum_{t\in T}\xi^2_t(\tau)\EE\left[\left(\ind[X_t\in A] - G_t(A)\right)^2\right] \label{eq:popo}\\
&\leq \frac{1}{4}\sum_{t\in T}\xi^2_t(\tau).
\]
For the second term, we first apply Jensen's inequality by noting that $\xi_t(\tau)\geq0$,  $\sum_{t\in T}\xi_t(\tau) = 1$,
\[
\EE\left[\left(\sum_{t\in T}\xi_t(\tau)G_t(A) - G_\tau(A)\right)^2\right] 
&\leq \sum_{t\in T}\xi_t(\tau)\EE\left[\left(G_t(A) - G_\tau(A)\right)^2\right].
\]
Since $0 \leq |G_t(A) - G_\tau(A)| \leq 1$, we have that $(G_t(A) - G_\tau(A))^2 \leq |G_t(A) - G_\tau(A)|$.
Finally by \cref{thm:dfrep}, we know that $\EE\left[\left|G_t(A) - G_\tau(A)\right|\right] \leq d(t,\tau)$. Hence
\[
\EE\left[\left(\sum_{t\in T}\xi_t(\tau)G_t(A) - G_\tau(A)\right)^2\right]
&\leq \sum_{t\in T}\xi_t(\tau)\EE\left[\left|G_t(A) - G_\tau(A)\right|\right]\\
&\leq \sum_{t\in T}\xi_t(\tau)d(t,\tau).
\]
We can combine the bounds on the first and second terms for each set $A_i$, $i\in\nats$, since $\sum_{i=1}^\infty c_i = 1$:
\[
\EE\left[\|\whG_{\tau} - G_{\tau}\|^2_\mcA\right] &\leq \frac{1}{4}\sum_{t\in T}\xi_t(\tau)^2 + \sum_{t\in T}\xi_t(\tau)d(t,\tau).\label{eq:prebound1}
\]
Before proceeding further with this bound by substituting the definition of $\xi_t(\tau)$, 
we will obtain a similar result for the tail bound. 
We add and subtract $\sum_{t\in T}\xi_t(\tau)G_t$ and use the triangle inequality:
\[
\Pr\left(\|\whG_{\tau} - G_{\tau}\|_\mcA > \delta\right) &\leq 
\Pr\left(\|\whG_{\tau} - \sum_{t\in T}\xi_t(\tau)G_t\|_\mcA +\| \sum_{t\in T}\xi_t(\tau)G_t - G_{\tau}\|_\mcA > \delta\right).
\]
By \cref{lem:rvsumbound}, we have
\[
\Pr\left(\|\whG_{\tau} - G_{\tau}\|_\mcA > \delta\right) &\leq 
\Pr\left(\left\|\whG_{\tau} - \sum_{t\in T}\xi_t(\tau)G_{t}\right\|_\mcA > \delta/2\right) +
\Pr\left(\left\|\sum_{t\in T}\xi_t(\tau)G_{t} - G_{\tau}\right\|_\mcA > \delta/2\right).
\]
For the first term in the sum, note that $\whG_\tau$ is a function of $X_T$,
which are conditionally independent given $G$.
Further note
that for each $t\in T$, the value of $\left\|\whG_{\tau} - \sum_{t\in T}\xi_t(\tau)G_{t}\right\|_\mcA$ can change by
at most $\xi_t(\tau)$ when varying the value of $X_t$.
Therefore by McDiarmid's inequality,
\[
\Pr\left(\left\|\whG_{\tau} - \sum_{t\in T}\xi_t(\tau)G_{t}\right\|_\mcA > \delta/2\right) 
&\leq \exp\left( - 2\frac{\left(\delta/2 - \EE\left\|\whG_{\tau} - \sum_{t\in T}\xi_t(\tau)G_{t}\right\|_\mcA \right)^2}{\sum_{t\in T}\xi_t(\tau)^2}\right),
\]
whenever $\delta/2 \geq \EE\left\|\whG_{\tau} - \sum_{t\in T}\xi_t(\tau)G_{t}\right\|_\mcA$.
Expanding the definition of the norm and using Jensen's inequality yields
\[
\EE\left\|\whG_{\tau} - \sum_{t\in T}\xi_t(\tau)G_{t}\right\|_\mcA &= \sum_{i=1}^\infty c_i \EE\left|\sum_{t\in T}\xi_t(\tau)(\ind[X_t\in A_i] - G_t(A_i))\right|\\
&\leq \sum_{i=1}^\infty c_i \sqrt{\EE\left(\sum_{t\in T}\xi_t(\tau)(\ind[X_t\in A_i] - G_t(A_i))\right)^2},
\]
at which point the same logic as in \cref{eq:popo} yields
\[
\EE\left\|\whG_{\tau} - \sum_{t\in T}\xi_t(\tau)G_{t}\right\|_\mcA
& \leq \frac{1}{2}\sqrt{\sum_{t\in T}\xi_t(\tau)^2},
\]
and hence for all $\delta \geq \sqrt{\sum_{t\in T}\xi_t(\tau)^2}$,
\[
\Pr\left(\left\|\whG_{\tau} - \sum_{t\in T}\xi_t(\tau)G_{t}\right\|_\mcA > \delta/2\right) 
&\leq \exp\left( -\frac{\left(\delta - \sqrt{\sum_{t\in T}\xi_t(\tau)^2}\right)^2}{2\sum_{t\in T}\xi_t(\tau)^2}\right).
\]
For the second term in the sum, we apply Markov's inequality to find that
\[
\Pr\left(\left\|\sum_{t\in T}\xi_t(\tau)G_{t} - G_{\tau}\right\|_\mcA > \delta/2\right)
&\leq 2\delta^{-1}\EE\left[\left\|\sum_{t\in T}\xi_t(\tau)G_{t} - G_{\tau}\right\|_\mcA\right]\\
&= 2\delta^{-1}\sum_{i=1}^\infty c_i\EE\left[\left|\sum_{t\in T}\xi_t(\tau)G_{t}(A_i) - G_{\tau}(A_i)\right|\right].
\]
Noting that $\sum_{t\in T}\xi_t(\tau) = 1$, we can apply Jensen's inequality,
\[
2\delta^{-1}\sum_{i=1}^\infty c_i\EE\left[\left|\sum_{t\in T}\xi_t(\tau)G_{t}(A_i) - G_{\tau}(A_i)\right|\right]
&= 2\delta^{-1}\sum_{i=1}^\infty c_i\EE\left[\left|\sum_{t\in T}\xi_t(\tau)\left(G_{t}(A_i) - G_{\tau}(A_i)\right)\right|\right]\\
&\leq 2\delta^{-1}\sum_{i=1}^\infty c_i\sum_{t\in T}\xi_t(\tau)\EE\left[\left|G_{t}(A_i) - G_{\tau}(A_i)\right|\right].
\]
Finally by local exchangeability and \cref{thm:dfrep},
\[
2\delta^{-1}\sum_{i=1}^\infty c_i\sum_{t\in T}\xi_t(\tau)\EE\left[\left|G_{t}(A_i) - G_{\tau}(A_i)\right|\right]
&\leq 2\delta^{-1}\sum_{i=1}^\infty c_i\sum_{t\in T}\xi_t(\tau)d(t,\tau)
= 2\delta^{-1}\sum_{t\in T}\xi_t(\tau)d(t,\tau).
\]
Combining the bounds for the first and second term and shifting $\delta$ yields, for all $\delta > 0$,
\[
\Pr\left(\|\whG_\tau - G_\tau\|_\mcA > \delta + \sqrt{\sum_{t\in T}\xi_t(\tau)^2}\right)
&\leq
\exp\left( - \frac{\delta^2}{2\sum_{t\in T}\xi_t(\tau)^2}\right)
+
2\frac{\sum_{t\in T}\xi_t(\tau)d(t,\tau)}{\delta + \sqrt{\sum_{t\in T}\xi_t(\tau)^2}}.\label{eq:prebound2}
\]
We now substitute the definition 
of $\xi_t(\tau) = \max\{0, 1/M_\tau + 2(\mu_\tau - d(t,\tau))\}$ into both results in \cref{eq:prebound1,eq:prebound2}.
First, note that (suppressing $(\tau)$ notation in the remainder 
of the proof for brevity), 
\[
\sum_{t\in T} \xi_t^2
&= \frac{1}{M_\tau} + 4M_\tau \sigma_\tau^2
&
&\text{and}
&
\sum_{t\in T}\xi_t d(t,\tau)
&= \mu_\tau -2M_\tau\sigma_\tau^2,
\]
where
\[
\mu_\tau &= \frac{1}{M_\tau}\sum_{m=1}^{M_\tau} d_m
&
&\text{and}
&
\sigma_\tau^2 &= \frac{1}{M_\tau}\sum_{m=1}^{M_\tau} d^2_m - \left(\frac{1}{M_\tau}\sum_{m=1}^{M_\tau}d_m\right)^2.
\]
Further, by \citet[Theorem 1]{Bhatia00} and the definition of $M_\tau$,
\[
\sigma_\tau^2 &\leq (\mu_\tau-d_1)(d_{M_\tau} - \mu_\tau) \leq \mu_\tau\left(\frac{1}{2M_\tau} +\mu_\tau - \mu_\tau\right) = \frac{\mu_\tau}{2M_\tau}.
\]
Therefore
\[
\EE\left[\|\whG_{\tau} - G_{\tau}\|^2_\mcA\right] &\leq 
\frac{1}{4M_\tau} + M_\tau\sigma_\tau^2 + \mu_\tau - 2M_\tau\sigma_\tau^2\\
&\leq \frac{1}{4M_\tau} + \mu_\tau,
\]
and
\[
\Pr\left(\|\whG_\tau - G_\tau\|_\mcA > \delta + \sqrt{\frac{1}{M_\tau}+4M_\tau\sigma_\tau^2}\right)
&\leq 
\exp\left( - \frac{M_\tau\delta^2}{2\left(1+4M_\tau^2\sigma_\tau^2\right)}\right)
+\frac{2\mu_\tau-4M_\tau\sigma_\tau^2}{\delta + \sqrt{\frac{1}{M_\tau}+4M_\tau\sigma_\tau^2}}\\
\Pr\left(\|\whG_\tau - G_\tau\|_\mcA > \delta + \sqrt{\frac{1}{M_\tau}+2\mu_\tau}\right)&\leq 
\exp\left( - \frac{M_\tau\delta^2}{2\left(1+2M_\tau\mu_\tau\right)}\right)
+\frac{2\mu_\tau}{\delta + \sqrt{\frac{1}{M_\tau}}}.
\]
Finally, because neither upper bound depends explicitly on $\mcA$, we can take the supremum.
To obtain the same results for $\Pr\left(X_\tau \in \cdot \given X_T\right)$, we
apply the same proof technique, noting that (1)
$\whG_\tau(A) = \EE\left[\whG_\tau(A) \given X_T\right]$ and (2)
 by the tower property 
and de Finetti result in \cref{thm:dfrep}, $\EE\left[\ind[X_{\tau}\in A] \given X_T\right] = \EE\left[G_\tau(A) \given X_T\right]$.
\eprfof

\bprfof{\cref{cor:empiricalmeasureasymp}}
For each $M \in [|T_n|]$, denote $\mu_M = \frac{1}{M}\sum_{m=1}^M d_m$. 
Note that for any $M < M_\tau$,
\[
&\left(\frac{1}{M}+2\mu_M - 2d_M\right) - \left(\frac{1}{M+1}+2\mu_{M+1} - 2d_{M+1}\right)\\
=&2(d_{M+1}-d_M) +
\left(\frac{1}{M}+2\mu_M \right) - \left(\frac{1}{M+1}+2\mu_{M+1}\right)\\
=&2(d_{M+1}-d_M) +
\frac{1}{M}\left(1+\sum_{m=1}^{M}2d_m\right) - \frac{1}{M+1}\left(1+\sum_{m=1}^{M+1}2d_m\right) \\
=& 2(d_{M+1}-d_M) +
\frac{1}{M+1}\left(\frac{1}{M} + \frac{1}{M}\sum_{m=1}^M 2d_m - 2d_{M+1}\right)\\
\geq& 2(d_{M+1}-d_M) +
\frac{1}{M+1}\left(2d_M - 2d_{M+1}\right)\\
\geq& \frac{2M}{M+1}(d_{M+1}-d_M).
\]
Therefore, for all $M< M_\tau$,
\[
\frac{1}{M}+2\mu_M - 2d_M \geq \frac{1}{M+1}+2\mu_{M+1}-2d_{M+1} + \frac{2M}{M+1}(d_{M+1}-d_M).
\]
We iterate this bound from $m=M$ to $m=M_\tau-1$ to find that for all $M < M_\tau$,
\[
\frac{1}{M}+2\mu_M - 2d_M &\geq \frac{1}{M_\tau} + 2\mu_\tau - 2d_{M_\tau} + \sum_{m=M}^{M_\tau - 1}\frac{2m}{m+1}(d_{m+1}-d_m).
\]
Finally, we rearrange this bound to obtain an upper bound on $\frac{1}{2M_\tau} + \mu_\tau$ for any $M < M_\tau$:
\[
\frac{1}{2M_\tau} + \mu_\tau &\leq \frac{1}{2M}+\mu_M - d_M + d_{M_\tau} - \sum_{m=M}^{M_\tau - 1}\frac{m}{m+1}(d_{m+1}-d_m)\\
&\leq \frac{1}{2M}+\mu_M - d_M + d_{M_\tau} - \frac{M}{M+1}\sum_{m=M}^{M_\tau - 1}(d_{m+1}-d_m)\\
&= \frac{1}{2M}+\mu_M - d_M + d_{M_\tau} - \frac{M}{M+1}(d_{M_\tau}-d_M)\\
&= \frac{1}{2M}+\mu_M + \frac{1}{M+1}(d_{M_\tau} - d_M)\\
&\leq \frac{1}{2M}+d_M + \frac{1}{M+1}.
\]
We also have that $M_\tau \to \infty$ as $n\to\infty$: by definition of $M_\tau$,
\[
\frac{1}{M_\tau+1} \leq \frac{1}{M_\tau+1} + 2\mu_{M_\tau+1} \leq 2d_{M_\tau+1},
\]
so if $\liminf_{n\to\infty} M_\tau = C < \infty$, then there would 
exist a subsequence such that $\frac{1}{C+1} \leq 2d_{C+1}$ for all $n$ sufficiently large.
But this is not possible, since for any fixed $C\in\nats$, $d_C \to 0$ as $T_n$ concentrates around $\tau$.
Therefore $M_\tau \to \infty$, so that for any $M\in\nats$,
\[
\limsup_{n\to\infty} \frac{1}{2M_\tau} + \mu_\tau &\leq \limsup_{n\to\infty} \frac{1}{2M}+d_M + \frac{1}{M+1}\\
&= \frac{1}{2M} + \frac{1}{M+1},
\]
and hence
\[
\limsup_{n\to\infty} \frac{1}{2M_\tau} + \mu_\tau = 0.
\]
\cref{thm:empiricalmeasurebds}  implies that both
$\EE\left[\|\whG_\tau - G_\tau\|^2_\mcA\right] \to 0$
and
$\EE\left[\|\whG_\tau - \Pr\left(X_\tau \in \cdot \given X_{T_n}\right)\|^2_\mcA\right] \to 0$
as $n\to\infty$.
By Markov's inequality, 
\[
\|\whG_\tau - G_\tau\|^2_\mcA \convp 0
\quad \text{and}\quad
\|\whG_\tau - \Pr\left(X_\tau \in \cdot \given X_{T_n}\right)\|^2_\mcA \convp 0, \qquad n\to\infty.\label{eq:tmpmarkov}
\]

Finally, note that \cref{eq:tmpmarkov} implies that any subsequence likewise satisfies $\|\cdot\|_\mcA \convp 0$, 
and hence any subsequence has a further subsequence such that $\|\cdot\|_\mcA \convas 0$.
Since $\mcA$ was arbitrary, \cref{lem:normAimpliesweakconv} asserts that we can choose $\mcA$ 
such that $\|\cdot\|_\mcA \to 0$ implies weak convergence, i.e., $\prokd{\cdot}{\cdot}\to 0$.
Thus any subsequence has a further subsequence 
that satisfies $\prokd{\cdot}{\cdot}\convas 0$. Hence $\prokd{\cdot}{\cdot} \convp 0$ by \cite[Theorem 2.3.2]{Durrett10}.
\eprfof

\bprfof{\cref{cor:diffusejumps}}
By Markov's inequality,
\[
\lim_{\eta \to 0}\sup_{t  \,:\, d(t,t_0)\leq \eta}\Pr\left(|G_t(A) - G_{t_0}(A)| > \epsilon\right)
&\leq \lim_{\eta \to 0}\sup_{t \,:\, d(t,t_0)\leq \eta} \epsilon^{-1} \EE\left|G_t(A) - G_{t_0}(A)\right|,
\]
and by \cref{thm:dfrep},
\[
&\leq \lim_{\eta \to 0}\sup_{t \,:\, d(t,t_0)\leq \eta} \epsilon^{-1} d(t,t_0)= 0.
\]
\eprfof

\bprfof{\cref{cor:stationarity}}
First, note that by assumption, the space $(d, \reals)$ is infinitely separable.
By local exchangeability and \cref{thm:dfrep},
for any $t, \Delta\in\reals$, finite subset $T\subset \reals$, and $A\in\Sigma$, \cref{thm:dfrep} implies that
\[
| \Pr(X_T \in A) - \Pr(X_{T+\Delta} \in A) |  \leq \sum_{t\in T} d(t, t+\Delta) = O(\Delta^{1+\gamma}), \quad \Delta \to 0, \label{eq:s1}\\
 \EE|G_t(A) - G_{t+\Delta}(A)| \leq  d(t, t+\Delta) = O(\Delta^{1+\gamma}), \quad \Delta\to 0, \label{eq:s2}
\]
where $T+\Delta$ denotes the translation of all covariates in $T$ by $\Delta$.
The Kolmogorov continuity theorem \cite[Theorem 3.23]{Kallenberg02} implies that
for all $\alpha\in(0,\gamma)$,  $(G_t(A))_{t\in\reals}$ has an $\alpha$-H\"older continuous modification.
Note that an $\alpha$-H\"older continuous function for $\alpha > 1$ is constant.

First, assume $\gamma > 1$. If we select $\alpha \in (1, \gamma)$, we have that for any $A\in\Sigma$,
$(G_t(A))_{t\in\reals}$ has a constant modification.  In other words, for all $t,t'\in\reals$, $A\in\Sigma$,
$\Pr\left(G_t(A) = G_{t'}(A)\right) = 1$. Since $\Sigma = \sigma(\mcA)$ 
for a countable algebra $\mcA$ \citep[Prop.~3.1, 3.3]{Preston08}, 
we have that $\Pr\left(\forall A \in \mcA, \,\, G_t(A) = G_{t'}(A)\right) = 1$, and hence by 
Carath\'eodory's extension theorem \citep[Theorem 2.5]{Kallenberg02}, $G_t$ and $G_{t'}$ are 
almost surely equal probability measures.
This implies that $G$ is a constant 
process (up to modification) and $X$ is exchangeable.

Next, suppose $\gamma \in (0, 1]$.  Then by \cref{eq:s1},
\[
\lim_{\Delta\to 0}\frac{| \Pr(X_T \in A) - \Pr(X_{T+\Delta} \in A) |}{\Delta}  \leq C\cdot \lim_{\Delta \to 0} \Delta^{\gamma} = 0,
\]
showing that $X$ is stationary. Next, since $X$ is stationary, for any $t, t' \in \reals$ and $A\in\Sigma$,
the mean of $G_t(A)$ satisfies
\[
&\EE\left[G_t(A)\right] = \Pr(X_t \in A) = \Pr(X_{t+\Delta} \in A) = \EE\left[G_{t+\Delta}(A)\right].
\]
Similarly, the autocovariance satisfies
\[
\MoveEqLeft{\EE\left[(G_t(A) - \EE G_t(A))(G_{t+\Delta}(A) - \EE G_{t+\Delta}(A))\right] }\\
&=\Pr(X_t\in A, X_{t+\Delta}\in A) - \Pr(X_t\in A)\Pr(X_{t+\Delta}\in A)\\
&=\Pr(X_{0}\in A, X_{\Delta}\in A) - \Pr(X_{0}\in A)\Pr(X_{\Delta}\in A)\\
&=\EE\left[(G_{0}(A) - \EE G_{0}(A))(G_{\Delta}(A) - \EE G_{\Delta}(A))\right].
\]
Hence $(G_t(A))_{t\in\reals}$ is weak-sense stationary.

Finally, consider the process $X_t = \ind(t \geq U)$ for $U\in \distUnif[0, 1]$,
which is locally exchangeable with $d(t,t') = \min(|t-t'|, 1)$ and hence $\gamma = 0$.
The underlying random measure process is specified by $G_t = \ind(t < U) \delta_{\{0\}} + \ind(t \geq U)\delta_{\{1\}}$
where $\delta_x$ is the Dirac measure at $x$; this has no sample-continuous modification. 
\eprfof

\bprfof{\cref{prop:approxsuff}}
Let $\mfT$ be the countable subset provided by infinite separability in \cref{assum:infsep}.
Let $(t_n)_{n=1}^\infty$ be any ordering 
of $\mfT \setminus T$, and $Y_N=\left(X_{t_N}, X_{t_{N+1}},\dots\right)$. 
Reverse martingale convergence implies that 
\[
\EE\left[h(X_T)|Y_N, \wtG\right] \convas \EE\left[h(X_T)|\mathcal{F}, \wtG\right]\eqas \EE\left[h(X_T)|G, \wtG\right] \,\, N\to\infty, \label{eq:convlimsuff}
\]
where $\mathcal{F}$ is the tail $\sigma$-algebra of $\{X_{t_i}\}_{i=1}^\infty$.
Defining $g(X_T)=\frac{1}{|\mcG|}\sum_{\pi\in \mcG}h(X_{\pi, T})$, we have that 
 $g(X_T)$ is invariant to $\mcG$ and thus $g(X_T)$ is $\sigma(\wtG,Y_N)$-measurable. Therefore
\[
\MoveEqLeft{\EE\left|\EE[h(X_T)|\wtG,Y_N]-g(X_T)\right|}\label{eq:convlimsuff2}\\
&=\EE\left[\frac{1}{|\mcG|}\left|\sum_{\pi\in \mcG}\EE[h(X_T)-h(X_{\pi, T})|\wtG,Y_N]\right|\right]\\
&\leq\EE\left[\frac{1}{|\mcG|}\sum_{\pi\in \mcG}\left|\EE\left[h(X_T)-h(X_{\pi, T})|\wtG,Y_N\right]\right|\right].
\]
By \cref{lem:conditionaldists}(1) and \cref{prop:finiteinfiniteequiv},
\[
\leq& 2\|h\|_\infty\EE\left[\frac{1}{|\mcG|}\sum_{\pi\in \mcG}d_{\mathrm{TV}}(X_T,X_{\pi, T})\right]\\
\leq& 2\|h\|_\infty\EE\left[\frac{1}{|\mcG|}\sum_{\pi\in \mcG}\sum_{t\in T}d(t,\pi(t))\right].
\]
Taking the limit as $N\to\infty$,  moving it into the expectation in \cref{eq:convlimsuff2} via dominated convergence, and using the limit from \cref{eq:convlimsuff} yields
\[
\EE\left|\EE[h(X_T)|\wtG,G]-g(X_T)\right|\leq 2\|h\|_\infty\EE\left[\frac{1}{|\mcG|}\sum_{\pi\in \mcG}\sum_{t\in T}d(t,\pi(t))\right].
\]
Identical reasoning to the above also shows that 
\[
\EE\left|\EE[h(X_T)|\wtG]-g(X_T)\right|\leq 2\|h\|_\infty\EE\left[\frac{1}{|\mcG|}\sum_{\pi\in \mcG}\sum_{t\in T}d(t,\pi(t))\right].
\]
Finally we add and subtract $g(X_T)$ in left hand side of \cref{eq:approxsuff}, apply the triangle inequality
with the above bounds, and note that the sum over $\pi$ is the expectation over a uniformly random permutation to obtain the result.
\eprfof

\bprfof{\cref{prop:exactrandomizationtest}}
We rewrite the probability as an expectation,
\[
&\Pr\left( \frac{1}{|\mcG|}\sum_{\pi \in \mcG} \ind\left[S(X_T) \leq S(X_{\pi,T})\right] \leq \alpha \right)\\
&= \EE\left[ \ind\left[\frac{1}{|\mcG|}\sum_{\pi \in \mcG} \ind\left[S(X_T) \leq S(X_{\pi,T}) \right] \leq \alpha \right] \right]\\
&= \frac{1}{|\mcG|}\sum_{\pi' \in \mcG} \EE\left[ \ind\left[\frac{1}{|\mcG|}\sum_{\pi \in \mcG} \ind\left[S(X_T) \leq S(X_{\pi,T})\right] \leq \alpha \right] \right].
\]
By local exchangeability, we can remap under any bijection $\pi' : T \to T$, so that
\[
&\leq \frac{1}{|\mcG|}\sum_{\pi' \in \mcG}\EE\left[\ind\left[\frac{1}{|\mcG|}\sum_{\pi \in \mcG} 
\ind\left[S(X_{\pi', T}) \leq S(X_{\pi\pi', T})\right] \leq \alpha \right] \right] + \sum_{t\in T}d(t, \pi'(t)).
\]
Finally, note that the outer indicator function tests whether $S(X_{\pi'T})$ is
strictly greater than $(1-\alpha) |\mcG|$ of the statistics across all
$\pi\in\mcG$. There can be at most $\alpha |\mcG|$ of such indicator functions, so
\[
&\leq \EE\left[\alpha \right] + \frac{1}{|\mcG|}\sum_{\pi' \in \mcG}\sum_{t\in T}d(t, \pi'(t))
= \alpha  + \frac{1}{|\mcG|}\sum_{\pi' \in \mcG}\sum_{t\in T}d(t, \pi'(t)).
\]
Rearranging the bound yields the result.
\eprfof

\section{Technical lemmata}\label{sec:techlemmas}

\bnlem\label{lem:conditionaldists}
Let $X, Y$ be bounded random variables in $[a,b]$ for some $a, b \in \reals$, $a\leq b$,
and $U, V$ be random elements in some probability space.
\benum
\item If $\|(X, U) - (Y, U)\|_{\mathrm{TV}} \leq \epsilon$, then
\[
\EE\left| \EE\left[X \given U\right] - \EE\left[Y \given U\right] \right|\leq (b-a)\epsilon.
\]
\item If $\|(X, U) - (X, V)\|_{\mathrm{TV}} \leq \epsilon$, then
for any 1-Lipschitz function $h : \reals\to \reals$,
\[
\left|\EE\left[h(\EE\left[X\given U\right]) - h(\EE\left[X\given V\right])\right]\right| \leq 3(b-a)\epsilon.
\]
\eenum
\enlem
\bprf
1. Denoting $Q \defined \ind\left[ \EE\left[X \given U\right] > \EE\left[Y \given U\right]\right]$,
\[
\EE\left| \EE\left[X \given U\right] - \EE\left[Y \given U\right] \right|
&= \EE\left[ \EE\left[X \given U\right]\left(2Q-1\right) - \EE\left[Y \given U\right] \left(2Q-1\right)\right].
\]
Using the fact that $Q$ is measurable with respect to $U$ and the tower property yields
\[
&= \EE\left[ X\left(2Q-1\right)\right] - \EE\left[Y \left(2Q-1\right)\right]\\
&= (b-a) \left(\EE\left[ \frac{X-a}{b-a}\left(2Q-1\right)\right] - \EE\left[\frac{Y-a}{b-a} \left(2Q-1\right)\right]\right).
\]
Since the difference is between the expectation of a function bounded in $[0, 1]$ evaluated at $(X,U)$ and at $(Y,U)$, the
assumed total variation bound provides the result.

2. First, note that $\sup_{x,y \in [a, b]} \left|h(x) - h(y)\right| \leq b-a$ by 1-Lipschitz continuity.
Then defining $A(U) \defined \EE\left[X \given U\right]$ and $B(V) \defined \EE\left[X \given V\right]$,
the triangle inequality yields
\[
\left|\EE\left[h(A(U))\right] - \EE\left[h(B(V))\right]\right|
&\leq\EE\left|A(U) - B(U)\right| + \left|\EE\left[h(B(V)) - h(B(U))\right]\right|.
\]
The right hand term is bounded by $(b-a)\epsilon$ by the assumed total variation bound and 1-Lipschitz continuity.
Defining $Q(u) = \ind\left[A(u) \geq B(u)\right]$, 
\[
\MoveEqLeft{\EE\left|A(U) - B(U)\right|}\\
 &=\EE\left[ A(U)(2Q(U)-1) - B(U)(2Q(U)-1)\right]\\
 &=(b-a)\EE\left[ \frac{A(U)-a}{b-a}(2Q(U)-1) - \frac{B(V)-a}{b-a}(2Q(V)-1) \right] \\
&+ (b-a)\EE\left[\frac{B(V)-a}{b-a}(2Q(V)-1) - \frac{B(U)-a}{b-a}(2Q(U)-1)\right].
\]
The first term in the expression can be bounded by $(b-a)\epsilon$ via substitution of the conditional expectation formulae for $A, B$, using the tower property, and controlling the difference
in expectations with the assumed total variation bound. The second term is again a difference in expectation of a bounded function under $U$ and $V$ with the same bound $(b-a)\epsilon$.
\eprf

\bnlem\label{lem:prodbd}
For any two sequences of real numbers $(a_i)_{i=1}^\infty$, $(b_i)_{i=1}^\infty$,
\[
\left|\prod_{i=1}^\infty a_i - \prod_{i=1}^\infty b_i\right| &\leq \sum_{i=1}^\infty |a_i-b_i|\left(\prod_{j=1}^{i-1}b_j\right)\left(\prod_{j=i+1}^\infty a_i\right).
\]
\enlem
\bprf
The proof follows by adding and subtracting $b_1\prod_{i=2}^\infty a_i$, then $b_1b_2\prod_{i=3}^\infty a_i$, etc., and then using the triangle inequality.
\eprf

\bnlem[\citep{Reiss81}] \label{lem:prodmsrbd}
For any two finite product probability
measures $\mu = \mu_1\times\dots \times \mu_N$ and
$\nu = \nu_1\times \dots \times \nu_N$,
\[
1-\exp\left(-\frac{1}{2}\sum_{n=1}^N d_{\mathrm{TV}}(\mu_n, \nu_n)^2\right)
\leq d_{\mathrm{TV}}(\mu,\nu)
\leq \sum_{n=1}^N d_{\mathrm{TV}}(\mu_n, \nu_n).
\]
\enlem

\bnlem\label{lem:rvsumbound}
For any two real-valued random variables $U,V$ and constants $a,b\in\reals$,
\[
\Pr\left(U+V > b\right) &\leq \Pr\left(U > a\right) + \Pr\left(V > b-a\right).
\]
\enlem
\bprf
\[
\Pr\left(U+V > b\right) &= \Pr\left(U+V > b\given U > a\right)\Pr\left(U>a\right) + \Pr\left(U+V > b \given U\leq a\right)\Pr\left(U\leq a\right)\\
&\leq \Pr\left(U>a\right) + \Pr\left(a+V > b \given U\leq a\right)\Pr\left(U\leq a\right)\\
&\leq \Pr\left(U>a\right) + \Pr\left(V > b-a, U\leq a\right)\\
&\leq \Pr\left(U>a\right) + \Pr\left(V > b-a\right).
\]
\eprf

\bnlem\label{lem:normAimpliesweakconv}
Let $(\mcX, \Sigma)$ be a standard Borel space.
There exists a countable collection of measurable subsets
$(A_i)_{i=1}^\infty$, $A_i \subseteq \mcX$
such that
for all $\mcA = \{c_i,A_i\}_{i=1}^\infty$, $c_i > 0$, $\sum_i c_i = 1$,
and
probability measures $\mu, (\mu_n)_{n=1}^\infty$, 
\[
\|\mu_n - \mu\|_\mcA \to 0 \implies \mu_n \convd \mu, \qquad n\to\infty,
\]
and for all $\mu$ such that each $A_i$ is a continuity set of $\mu$,
\[
\mu_n \convd \mu \implies \|\mu_n - \mu\|_\mcA \to 0, \qquad n\to\infty.
\]
\enlem
\bprf
Since $(\mcX, \Sigma)$ is a standard Borel space,
we know that $\Sigma$ is generated by a topology with a countable base
$(B_i)_{i=1}^\infty$. Any open set $U \subseteq \mcX$ can be expressed
as a countable union of these sets.
Consider the collection of all possible unions $\mcB_n$ of $\{B_1, \dots, B_n\}$,
and construct a countable sequence of sets $(A_i)_{i=1}^\infty$ by ordering $\mcB_1$, then $\mcB_2$, and so on.
Then for any open set $U\subseteq \mcX$, there exists a subsequence 
$(U_k)_{k=1}^\infty$ of $(A_i)_{i=1}^\infty$ such that $U_k \uparrow U$.

Assume $\|\mu_n - \mu\|_\mcA \to 0$;
then for any open set $U$ and $k\in\nats$,
\[
\liminf_{n\to\infty}\mu_n(U) \geq \liminf_{n\to\infty} \mu_n(U_k).
\]
But $\|\mu_n - \mu\|_{\mcA} \to 0$ if and only if $\forall i\in\nats$, $\mu_n(A_i) \to \mu(A_i)$.
Hence
\[
\liminf_{n\to\infty} \mu_n(U_k) = \mu(U_k).
\]
Since this holds for all $k\in\nats$ and $U_k\uparrow U$, by the continuity of measures,
\[
\liminf_{n\to\infty}\mu_n(U) \geq \mu(U).
\]
Hence $\mu_n \convd \mu$. If each $A_i$ is a continuity set of $\mu$, then
$\mu_n \convd \mu$ implies that $|\mu_n(A_i) - \mu(A_i)| \to 0$ for each $i$,
which then implies $\|\mu_n - \mu\|_\mcA \to 0$.
\eprf

\section{Additional Examples}\label{sec:examples_add} 
In this section, we show that many popular covariate-dependent models 
from Bayesian nonparametrics exhibit local exchangeability.

\subsection{Dependent Dirichlet process mixtures}
In a typical mixture model setting, we have observations generated via
\[
X_n&\distiid \sum_{k=1}^\infty w_k F(\cdot; \theta_k),  \qquad n\in\nats,
\]
where $(w_k)_{k=1}^\infty$ are the mixture weights satisfying $w_k \geq 0$, $\sum_k w_k = 1$; 
$(\theta_k)_{k=1}^\infty$ are the component parameters;
$F(\cdot ; \theta)$ is the mixture component likelihood;
and $(X_n)_{n=1}^\infty$ are the observations. A popular nonparametric
prior for the weights and component parameters is the Dirichlet process \citep{Ferguson73}, 
defined by \citep{Sethuraman94}
\[
\theta_k\distiid H, \qquad v_k\distiid \distBeta(1, \alpha), \qquad w_k = v_k\prod_{i=1}^{k-1}(1-v_i), \qquad k&\in\nats,
\]
for some distribution $H$. When the observations come with additional covariate information,
the \emph{dependent} Dirichlet process mixture model \citep{MacEachern99,MacEachern00} may be used to
  capture similarities between related mixture population data. Here, observations are generated via
\[
X_{x,n}&\distind \sum_{k=1}^\infty w_{x,k} F(\cdot; \theta_{x,k}), \qquad n\in\nats, \,\,x\in\reals,
\]
where the component parameters $\theta_{x,k}$ and stick variables $v_{x,k}$ are now \iid stochastic 
processes on $\reals$, 
and  $w_{x,k} = v_{x,k}\prod_{i=1}^{k-1}(1-v_{x,i})$. 
The marginal distributions of $\theta_{x,k}$ and $v_{x,k}$ at $x\in\reals$ are $H$ and $\distBeta(1,\alpha)$, respectively.
Thus, the dependent Dirichlet process is marginally a Dirichlet 
process for each covariate value, but can exhibit a wide range of dependencies across covariates. In this setting, we have $\mcT = \reals\times \nats$
and strong canonical premetric 
\[
d_{sc}(t,t') &= \EE\left[d_{\mathrm{TV}}\left(\sum_{k=1}^\infty w_{x,k} F(\cdot; \theta_{x,k}), \sum_{k=1}^\infty w_{x',k} F(\cdot; \theta_{x',k})\right)\right]\\
&= \frac{1}{2}\EE\left[\int \left|\sum_{k=1}^\infty w_{x,k} F(y; \theta_{x,k}) - \sum_{k=1}^\infty w_{x',k} F(y; \theta_{x',k})\right| \dee y\right],
\]
where $t = (x,n)$ and $t'=(x',n')$.
We add and subtract $\sum_{k=1}^\infty w_{x',k} F(\cdot; \theta_{x,k})$ and
apply the triangle inequality to find that 
\[
d_{sc}(t,t') &\leq 
\EE\left[d_{\mathrm{TV}}(F(\cdot; \theta_{x,1}), F(\cdot; \theta_{x',1}))\right]+
\sum_{k=1}^\infty \EE\left|w_{x,k}-w_{x',k}\right|.
\]
Since $w_{x,k}$ is a product of independent random variables, \cref{lem:prodbd} yields
\[
d_{sc}(t,t') &\leq 
\EE\left[d_{\mathrm{TV}}(F(\cdot; \theta_{x,1}), F(\cdot; \theta_{x',1}))\right]+\\
&\EE\left[\left|v_{x,1}-v_{x',1}\right|\right]\sum_{k=1}^\infty 
\left(\left(\frac{\alpha}{\alpha+1}\right)^{k-1} +\frac{k-1}{1+\alpha}\left(\frac{\alpha}{\alpha+1}\right)^{k-2}\right).
\]
The infinite sum converges to some $0<C<\infty$, and so
\[
d_{sc}(t,t') &\leq \EE\left[d_{\mathrm{TV}}(F(\cdot; \theta_{x,1}), F(\cdot; \theta_{x',1}))\right]+ C\EE\left|v_{x,1}-v_{x',1}\right|.
\]
Therefore, if the stochastic processes for the parameters and stick variables are
both smooth enough such that 
\[
\max \left\{ \EE\left[d_{\mathrm{TV}}(F(\cdot; \theta_{x,1}), F(\cdot; \theta_{x',1}))\right], \EE\left|v_{x,1}-v_{x',1}\right|\right\} \leq (1+C)^{-1}\td(t,t'),
\]
for some premetric $\td : \mcT \times \mcT \to \reals_+$, then $X$ is locally exchangeable
with respect to $\min(1, \td)$.
Many dependent processes (e.g., \citep{Foti15}) similar to the dependent Dirichlet process (and kernel beta process below) can be shown to exhibit local exchangeability using similar techniques.

\subsection{Kernel beta processes} Another example of a model exhibiting local exchangeability from the Bayesian nonparametrics
literature is the kernel beta process latent feature model \citep{Ren11}. In a typical nonparametric latent
feature modelling setting, we have observations generated via
\[
X_n = F\left(\cdot; Z_n\right), \qquad Z_n \distind \distBeP(\sum_{k=1}^\infty w_k \delta_{\theta_k}),
\] 
where $(w_k)_{k=1}^\infty$ are the feature frequencies satisfying $w_k \in [0, 1]$, $\sum_{k=1}^\infty w_k < \infty$;
$(\theta_k)_{k=1}^\infty$ are the feature parameters;
$\distBeP$ is the Bernoulli process that sets $Z_n(\{\theta_k\}) = 1$ with probability $w_k$ and $0$ 
otherwise independently across $k\in\nats$;
and $F$ is the likelihood for each observation. A popular nonparametric prior for the 
weights and feature parameters is the beta process \citep{Hjort90}, defined by 
\[
(\theta_k, w_k)_{k=1}^\infty \dist \distPP(\gamma H(\dee \theta) c(\theta)w^{-1}(1-w)^{c(\theta)-1}\dee w),
\]
where $\distPP$ is a Poisson point process parametrized by its mean measure, $c$ is some  positive function, $H$ is a probability distribution, and $\gamma > 0$. 
When the observations come with covariate information, the  
\emph{kernel} beta process \citep{Ren11} may be used to capture similarities in the latent features of related populations.
In particular, we replace $Z_n$ with
\[
Z_{x,n} \distind \distBeP\left(\sum_{k=1}^\infty \kappa(x, x_k; \psi_k)w_k\delta_{\theta_k}\right),
\] 
where $\kappa(x, x_k; \psi_k)$ is a kernel function with range in $[0, 1]$ centered at $x_k$ with parameters $\psi_k$, and
\[
(x_k, \psi_k, \theta_k, w_k)_{k=1}^\infty \dist \distPP(Q(\dee x)R(\dee \psi)\gamma H(\dee \theta) c(\theta)w^{-1}(1-w)^{c(\theta)-1}\dee w),
\]
where $Q$ and $R$ are probability distributions. In other words, the kernel
beta process endows each atom with \iid covariates $x_k$ and parameters
$\psi_k$, and makes the likelihood that an observation with covariate $x$ 
selects a feature with covariate $x_k$ depend on both $x$ and $x_k$. Taking $\reals$ to be the
space of covariates
for simplicity, again we have $\mcT = \reals\times\nats$ and (marginalizing $Z_{x,n}$) strong canonical premetric 
\[
d_{sc}(t,t') = \EE\left[d_{\mathrm{TV}}\left(\EE\left[F(\cdot; Z_{x,n})\right], \EE\left[F(\cdot; Z_{x',n'})\right]\right)\right],
\]
where $t = (x,n)$ and $t'=(x',n')$.
Suppose $F$ is $\gamma$-H\"older continuous in total variation for $0<\gamma \leq 1$, $C\geq 0$ in the sense that
\[
d_{\mathrm{TV}}\left(\EE\left[F(\cdot; Z)\right], \EE\left[F(\cdot ; Z')\right]\right) &\leq C\left(\sum_{k=1}^\infty \|\theta_k\|\left|p_k-p'_k\right| \right)^\gamma
\]
for any collection of points $\{\theta_k\}_{k=1}^\infty$, where $Z(\{\theta_k\}) = 1$, $Z'(\{\theta_k\}) = 1$ independently with probability $p_k$ and $p'_k$, respectively, and both assign 0 mass
to all other sets.
Then 
\[
d_{sc}(t,t') 
\leq C\EE\left[\sum_{k=1}^\infty |\kappa(x, x_k; \psi_k)-\kappa(x', x_k; \psi_k)|w_k\left\|\theta_k\right\|\right]^\gamma.
\]
Finally, if the kernel $\kappa$ is $\alpha$-H\"older continuous with constant $C'(\psi)$ depending on $\psi$,
 the independence of $\theta_k$, $w_k$, and $\psi_k$ may be used to show that
\[
d_{sc}(t,t') 
&\leq C\EE\left[\sum_{k=1}^\infty C'(\psi_k)|x-x'|^{\alpha}w_k\left\|\theta_k\right\|\right]^\gamma\\
&= C\left(\EE\left[C'(\psi_1)\right]|x-x'|^{\alpha}\EE\left[\left\|\theta_1\right\|\right]\EE\left[\sum_{k=1}^\infty w_k\right]\right)^\gamma.
\]
Therefore the observations are locally exchangeable with
$d(t, t') = \min\left(1, C''|x-x'|^{\alpha\gamma}\right)$ and $C''$ collects the product of constants from the previous expression.

\subsection{Dynamic topic model}
The dynamic topic model \citep{Blei06,Wang08} is a model for text data 
that extends latent Dirichlet allocation \citep{Blei03} to incorporate timestamp covariate information.
In a continuous version of the model, observations are generated via
\[
 D_{n,x} \dist \distMulti(W, \sum_{k=1}^K\theta_{x,k}\pi_V(\beta_{x,k})), \qquad \theta_x \dist \distDir(\pi_K(\alpha_x)), \qquad W \dist \distPoiss(\mu),
\]
where $x\in\reals$ represents timestamps, $\alpha_x\in\reals^K$ is a vector of
$K$ independent Wiener processes representing the popularity of $K$ topics at
time $x$, $\beta_{x,k}\in\reals^V$ is a vector of $V$ independent Wiener
processes representing the word frequencies for vocabulary of size $V$ in topic
$k$, $\pi_J$ is any $L$-Lipschitz mapping from $\reals^J$ to the probability
simplex $\pi_J : \reals^{J} \to \Delta^{J-1}$ for any $J\in\nats$, $\mu$ is the
mean number of words per document, $D_{n,x}\in\nats^V$ is the vector of
counts of each vocabulary word in the $n^\text{th}$ document observed at time
$x$, and $W$ is the number of words in each document, taken to be the same across
documents for simplicity.  Here the covariate space is $\mcT = \reals\times\nats$, and the
observations are count vectors in $\nats^V$ where $V$ is the vocabulary size.
In this setting, the strong canonical premetric is
\[
d_{sc}(t,t') &= \EE\left[d_{\mathrm{TV}}\left(
\distMulti(W, \sum_{k=1}^K\theta_{x,k}\pi_V(\beta_{x,k})),\distMulti(W, \sum_{k=1}^K\theta_{x',k}\pi_V(\beta_{x',k}))\right)\right],
\]
where $t=(x,n)$ and $t'=(x',n')$.
But since multinomial variables are a function (in particular, a sum) of independent categorical random variables, \cref{lem:prodmsrbd} yields the bound
\[
 d_{sc}(t,t')&\leq \EE\left[Wd_{\mathrm{TV}}\left(
\distCat(\sum_{k=1}^K\theta_{x,k}\pi_V(\beta_{x,k})),\distCat(\sum_{k=1}^K\theta_{x',k}\pi_V(\beta_{x',k}))\right)\right].
\]
We evaluate the total variation between two categorical distributions and apply the triangle inequality to find that
\[
 d_{sc}(t,t')&\leq \EE\left[\frac{W}{2}\sum_{v=1}^V\left|\sum_{k=1}^K\theta_{x,k}\pi_V(\beta_{x,k})_v - \sum_{k=1}^K\theta_{x',k}\pi_V(\beta_{x',k})_v\right|\right]\\
&\leq \frac{\mu}{2}\sum_{v=1}^V\sum_{k=1}^K\EE\left[\left|\theta_{x,k}-\theta_{x',k}\right|\pi_V(\beta_{x,k})_v + \theta_{x',k}\left|\pi_V(\beta_{x,k})_v - \pi_V(\beta_{x',k})_v\right|\right].
\]
Since $\sum_{v=1}^V \pi_V(\beta_{x,k})_v = \sum_{k=1}^K \theta_{x',k} = 1$, the components of $\theta_{x,k}$ and $\beta_{x,k}$ are \iid across $k$,
and $\pi_V$ is $L$-Lipschitz,
\[
 d_{sc}(t,t')&\leq \frac{\mu}{2}\left(KL \EE\left|\alpha_{x,1}-\alpha_{x',1}\right| + VL\EE\left|\beta_{x,1,1} - \beta_{x',1,1}\right|\right]\\
&\leq \frac{\mu L \left(K+V\right)}{2}\sqrt{|x-x'|},
\]
where the last line follows by Jensen's inequality.
Therefore the observations are locally exchangeable with $d(t, t') = \min\left(1, \frac{1}{2}\mu L \left(K+V\right)\sqrt{|x-x'|}\right)$.

\small
\bibliographystyle{imsart-nameyear}
\bibliography{sources}

\end{document}